\documentclass[final,5p,times]{elsarticle}

\usepackage{graphicx}
\newcommand{\dt}{\mathrm{d}t}

\newcommand{\dd}{\mathrm{d}}
\usepackage{dsfont}
\usepackage{lineno}
\usepackage[hidelinks]{hyperref}
\usepackage{amsmath}
\usepackage{amssymb}
\usepackage{wrapfig}
\usepackage{comment}
\usepackage{subfigure}
\usepackage{floatflt}
\usepackage{dsfont}
\usepackage{overpic}
\usepackage{tikz}
\usepackage[siunitx]{circuitikz}
\usepackage{pgfplots}
\pgfplotsset{width=7cm,compat=1.3}
\usepackage{mwe}
\usepackage{calc}
\usepackage{xcolor}
\usepackage{footnote}
\makesavenoteenv{tabular}
\makesavenoteenv{table}
\usepackage{pgfplots} 
\usepackage{url}
\pgfplotsset{compat=newest} 
\pgfplotsset{plot coordinates/math parser=false} 
\newlength\figureheight 
\newlength\figurewidth 

\biboptions{sort&compress}

\modulolinenumbers[1]

\makeatletter
\def\ps@pprintTitle{%
  \let\@oddhead\@empty
  \let\@evenhead\@empty
  \let\@oddfoot\@empty
  \let\@evenfoot\@oddfoot
}
\makeatother

\begin{document}

\begin{frontmatter}

\title{Analysis and simulation of a modified cardiac cell model gives \\ accurate predictions of the dynamics of the original one}

\author[AE1]{Andr\'e H. Erhardt}
\address[AE1]{Weierstrass Institute for Applied Analysis and Stochastics, Mohrenstraße 39, 10117 Berlin, Germany}
\ead{andre.erhardt@wias-berlin.de}
\author[SS]{Susanne Solem}
\address[SS]{Department of Mathematics, Norwegian University of Life Sciences, Norway}


\begin{abstract}
The 19-dimensional TP06 cardiac muscle cell model is reduced to a 17-dimensional version, which satisfies the required conditions for performing an analysis of its dynamics by means of bifurcation theory. The reformulated model is shown to be a good approximation of the original one. As a consequence, one can extract fairly precise predictions of the behaviour of the original model from the bifurcation analysis of the modified model. Thus, the findings of bifurcations linked to complex dynamics in the modified model --- like early afterdepolarisations (EADs), which can be precursors to cardiac death --- predicts the occurrence of the same dynamics in the original model. It is shown that bifurcations linked to EADs in the modified model accurately predicts EADs in the original model at the single cell level. Finally, these bifurcations are linked to cardiac death at the tissue level by example.   
 
 
\end{abstract}

\begin{keyword} 
nonlinear dynamics \sep reaction-diffusion system \sep cardiac muscle cell \sep mathematical modelling
\MSC[2010] 37G15 \sep 37N25 \sep 35Q92 \sep 65P30 \sep 92B05
\end{keyword}

\end{frontmatter}


\section{Introduction} 
Mathematical modelling and simulations are important to investigate and analyse phenomena in life science. Two advantages of a detailed mathematical model is the possibility to use mathematical theory in order to get an in-depth understanding of the underlying mechanisms of certain dynamics, and to develop numerical experiments instead of real experiments. For instance, the mathematical models may be used for drug testing and to study drug response or its properties~\cite{drugs1,drugs2,drugs3,drug4,drug5,Rodriguez_drug}. 
To derive a suitable mathematical model describing complex behaviour in biology like a cardiac action potential, based on experimental data, displaying all facets and dynamics of a cardiac muscle cell is a tough challenge. In addition to normal action potentials of a cardiac muscle cell, certain kinds of cardiac arrhythmia or chaos may occur~\cite{DELANGE2012365,Roden,VN,Vandersickel1,Vandersickel2}. This includes specific types of abnormal heart rhythms, which can lead to sudden cardiac death. It is therefore highly interesting and important to understand the complex behaviour and mechanism of such biological phenomena. A good mathematical model in combination with a meaningful analysis helps to decode these phenomena, the most prominent being early afterdepolarisations (EADs). However, due to the above mentioned multifaceted nature of cardiac cells, cardiac cell models can indeed be quite complex \cite{modeloverview,FINK20112}. This poses challenges both from a mathematical and computational point of view. 


The (endocardial) TP06 cardiac cell model \cite{TP06}, which is based on the TNNP04 model~\cite{TNNP04}, is by now a well-established and much studied cardiac cell model \cite{VN,Vandersickel1,Vandersickel2,modeloverview,FINK20112,Qu2021}. Despite that, a rigorous mathematical analysis of all its dynamics is still lacking. Performing such an analysis on the model is hard due to two main things: the model is 19-dimensional and it is not smooth. 

In this paper we modestly alter the non-smooth TP06 model in order to generate a smooth 17-dimensional version for which we can apply bifurcation theory. We show that the bifurcation theory aids in identifying complex dynamics linked to cardiac arrhythmia of the modified model, such as EADs. Through several examples, both at the single cell level and at the tissue level, we show that the bifurcation analysis of the modified model predicts similar, or even identical, dynamics of the original TP06 model. Thus, a bifurcation analysis of the modified model can help in decoding complex dynamics of the original model, which in turn can further the understanding of the heart.

Finally, we would like to emphasise that the approach in this paper can be utilised for other similar models stemming from the life sciences.





\section{On modelling cardiac action potentials and related issues}
\begin{figure}[h]
\centering
\includegraphics[width=0.75\columnwidth]{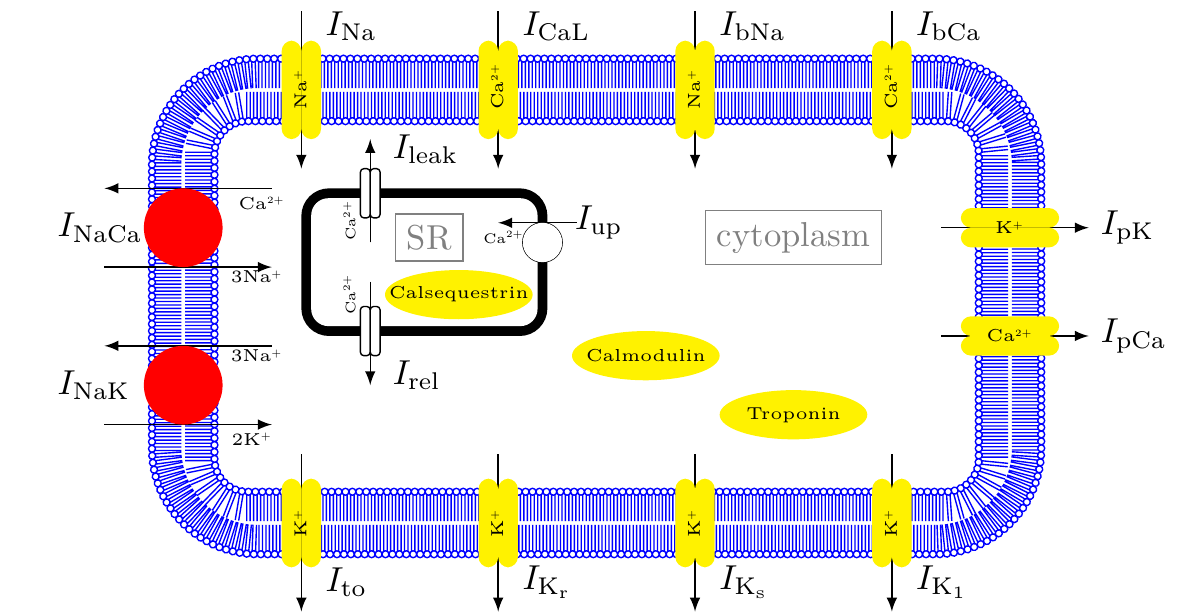}
\caption{Scheme of a cardiac myocyte: SR denotes the sarcoplasmic reticulum, $I_\text{NaCa}=\mathrm{Na}^+/\mathrm{Ca}^{2+}$ a exchanger current \&  $I_\text{NaK}=\mathrm{Na}^+/\mathrm{K}^+$ a pump current~\cite{TNNP04}.}\label{fig:cell}
\end{figure}

The mathematical modelling of action potentials (APs) of excitable biological cells like neurons and cardiac muscle cells has its origin in the famous and pioneering Hodgkin--Huxley model
~\cite{HH}. Here, an approach was established that can be utilised to model APs of excitable biological cells by a system of ordinary differential equations (ODEs). In the last decades there has been an immense development in the modelling of cardiac muscle cells.  
These conductance-based models represent a minimal biophysical interpretation of an excitable biological cell in which current flow across the membrane is due to charging of the membrane capacitance and movement of ions across ion channels, which are selective for particular ionic species, cf. Figure~\ref{fig:cell}. An initial stimulus activates the ion channels as soon as a certain threshold potential is reached. Then, these ion channels break open and/or up allowing an ion current flow, which changes the membrane potential. 
This electrophysiological behaviour can be described by the following ODE using the Hodgkin--Huxley formalism: 
\begin{align}\label{ode}
	C_m\frac{\dd V}{\dt}=-I_\mathrm{ion}+I_\mathrm{stim}.
\end{align}
Here $V$ denotes the voltage (in $mV$) and $t$ the time (in $ms$), while $I_\mathrm{ion}$ is the sum of all transmembrane ionic currents. $I_\text{stim}=52~pA/pF$ (applied for $2~ms$) represents the externally applied initial stimulus and $C_m=1~\mu F/cm^2$ the membrane capacitance.



 The main issue with several of the cardiac models of the form \eqref{ode}, from a mathematical point of view, is the requirement of a sufficiently smooth and regular ODE system in order to utilise theory such as bifurcation analysis. Bifurcation theory has already been widely utilised to investigate the dynamics of cardiac muscle cell models, see~\cite{AE_control,AE_MMOs,Kurata,OTTE2016265,Tran,Tsumoto2017,herrero2020reduced,Rose,DIEKMAN2021319}. However, many of the existing detailed models are not smooth enough for the theory to be applied, see for example~\citep{Priebe}. A second problem is that these systems may not exhibit a (suitable) equilibrium, i.e. the long-time behaviour ($t\to\infty$) of the system does not converge into a stable steady state.

In addition, very detailed and high-dimensional models also cause numerical issues, as they might be quite sensitive and computationally expensive. This makes them harder to investigate. In general, lower-dimensional models are easier and faster to analyse and to simulate. However, important information can be lost. The desired aim of cardiac cell modelling is to study a most detailed cardiac muscle cell model, which simultaneously is mathematically and computationally feasible. To this end, we will focus on the 19-dimensional TP06 model~\cite{TP06}, a cardiac model of the form \eqref{ode} exhibiting the above mentioned issues: lack of smoothness, high-dimensionality, and lack of a computationally stable equilibrium. We will modestly reformulate the model in order to derive a good approximation and to be able to apply numerical bifurcation analysis using the \texttt{CL\_MATCONT}, \textit{a continuation toolbox for MATLAB} \cite{Dhooge,Dhooge1,Govaerts}. For this modified TP06 model symbolic derivatives up to the fifth order can be derived, which is desirable as it allows for discovering more complex bifurcations using the algorithm.   

\begin{figure}[h]
\centering
\begin{overpic}[width=0.8\columnwidth]{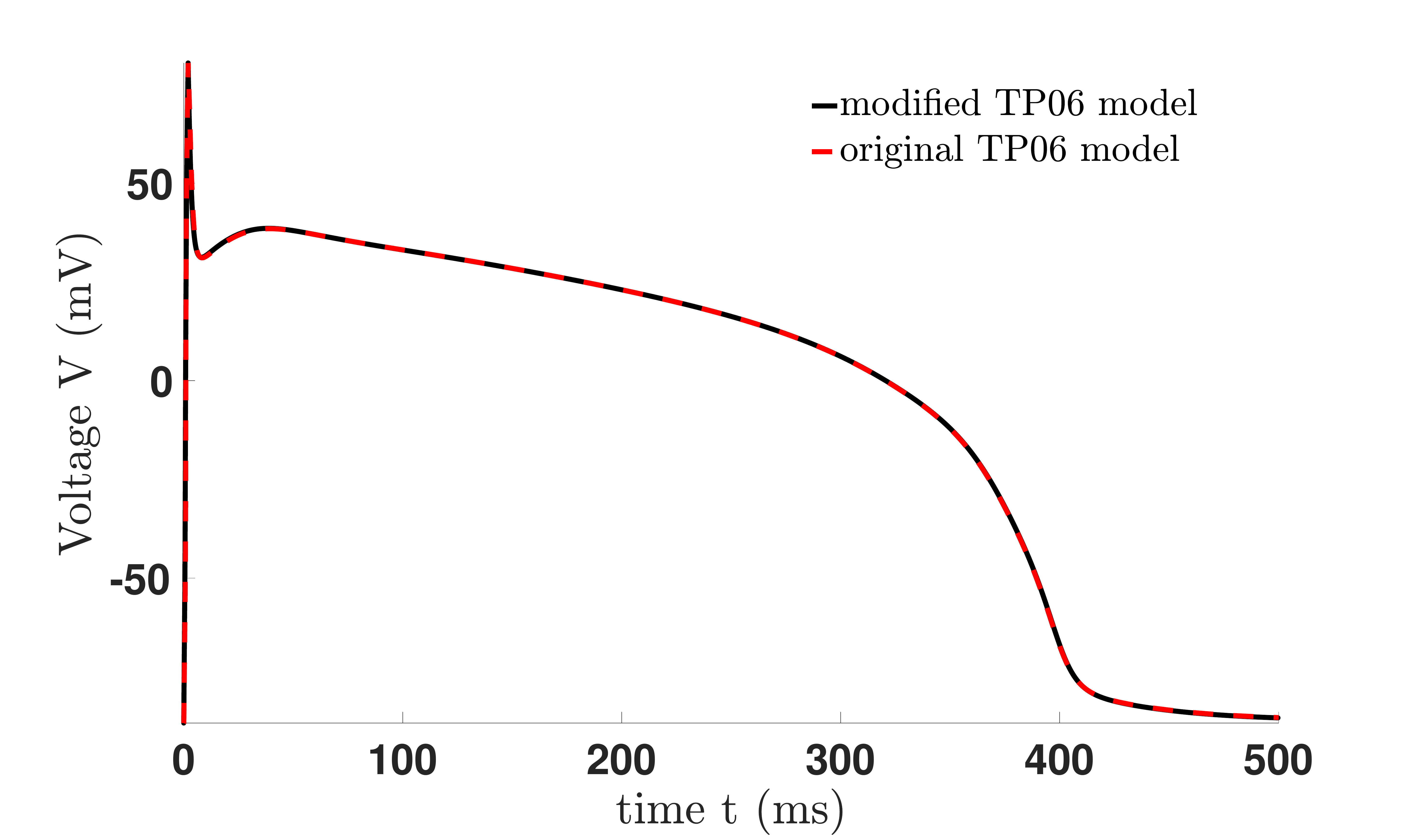}
\put(0,50){(a)}
\end{overpic}
\begin{overpic}[width=0.8\columnwidth]{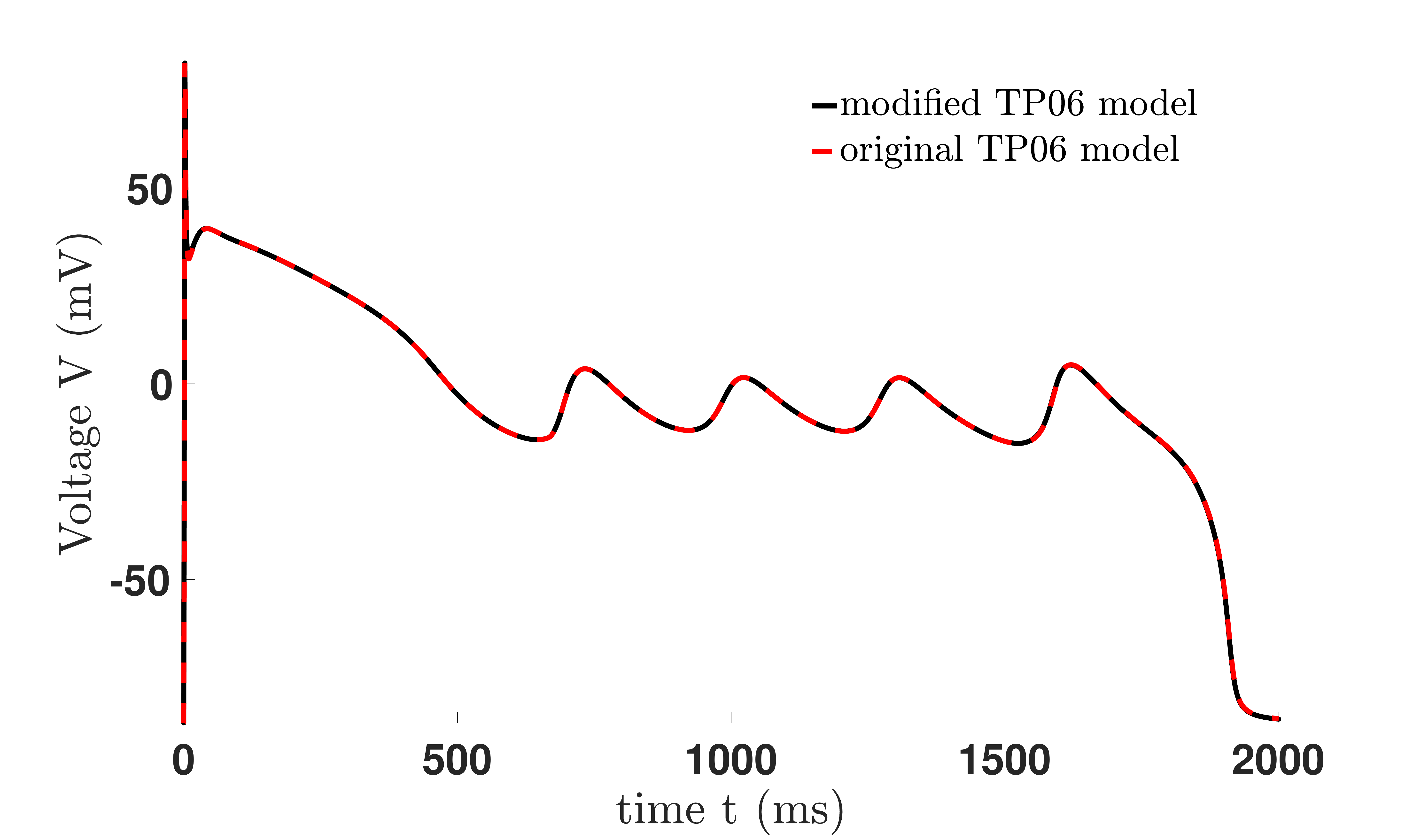}
\put(0,50){(b)}
\end{overpic}
\caption{Comparison of the trajectories of the original TP06 model with the modified one. (a) normal action potential. (b) early afterdepolarisation. Here, we use $G_\text{CaL}=0.000175~cm^3/(\mu Fs)$, cf.~\cite{TNNP04}, to adjust the AP duration.}\label{fig:comparison}
\end{figure}

Although it is a modified TP06 model which is investigated using bifurcation theory in this paper, we will present several examples validating that the analysis presented also provides fairly precise predictions of the dynamics of the TP06 model itself. As a first example, the two models produce indistinguishable output in the two plots in Figure~\ref{fig:comparison}. Figure~\ref{fig:comparison}(b) illustrates one example of an EAD. To be more precise, EADs are pathological voltage oscillations during the plateau or the repolarisation phase of the cardiac action potential and considered as potential precursors to cardiac arrhythmia.

\vspace{-0.5cm}
\section{A reformulation of the TP06 model}
The TP06 model contains several different ion currents, ion pump, ion exchanger and background currents:
\begin{align*}
I_\text{ion}=I_\text{K1}+I_\text{to}+I_\text{Kr}+I_\text{Ks}+I_\text{CaL}+I_\text{NaK}+I_\text{Na}
\\
+I_\text{bNa}+I_\text{NaCa}+I_\text{bCa}+I_\text{pK}+I_\text{pCa},
\end{align*}
cf. \cite{TNNP04,TP06} and Figure~\ref{fig:cell}. These currents are depending on individual ionic conductances $G_\text{current}$ and Nernst potentials $E_\text{current}$. Moreover, they may depend on gating variables, which are important for the activation and inactivation of the ion currents. For a full description of the variables involved, see \cite{TP06}. The parameters used in this paper are the same as those in~\cite{TP06}, with the exception of the $G_\text{CaL}$ value, which is set to $G_\text{CaL} = 0.000175~cm^3/(\mu Fs)$. 

The first smoothness issue is in the modelling of the sodium current in \cite{TNNP04,TP06}, i.e.
$$
I_\text{Na}=G_\text{Na}m^3hj(V-E_\text{Na}),
$$
where $G_\text{Na}$ denotes the ionic conductance and $E_\text{Na}$ the Nernst potential, while the different gating variables $m$, $h$ and $j$ satisfy 
$$
	\frac{\dd g}{\dt}=a_g(1-g)-b_gg=a_g-(a_g+b_g)g=\frac{g_\infty-g}{\tau_g},
$$
where $g$ represents the gating variables $m$, $h$ and $j$, while $g_\infty:=g_\infty(V)=a_g\cdot(a_g+b_g)^{-1}$ denotes the equilibrium of the gating variable $g$ and $\tau_g:=\tau_g(V)=(a_g+b_g)^{-1}$ its time scale. To be more precise, in the modelling of the gating variables $h$ and $j$, the voltage dependent functions $a_h$, $b_h$, $a_j$ and $b_j$ are not continuous, cf. \cite{TNNP04}. As the equilibrium of $h$ and $j$ are equal, we use the idea from~\cite{Bernus} and reformulate $h$ and $j$ to one new gating variable $v$, and modify the sodium current to
$$
I_\text{Na}=G_\text{Na}m^3v^2(V-E_\text{Na}).
$$
Now, $v$
satisfies the differential equation from above with $v_\infty=h_\infty=j_\infty$ and the time relaxation constant is given by
$$
\tau_v=0.25+\frac{2.24\cdot v_\infty}{(1-\tanh(6.468+0.07\cdot V))}.
$$
As can be seen in the figures throughout the paper, this modification does not change the behaviour of the system dramatically. However, it is indeed a reformulation of the original TP06 model and the behaviours are not expected to be identical in all cases. Notably is also, regarding the bifurcation analysis, that we start with a steady state solution or equilibrium of the investigated system, i.e. $v_\infty=h_\infty=j_\infty$ is naturally satisfied.

A further problem is evoked by the modelling of the intracellular potassium ion concentration $[K]_i$. The modelling is the standard one, i.e.
\begin{align*}
\frac{\dd[K]_i}{\dt}=-\frac{I_\text{K1}+I_\text{to}+I_\text{Kr}+I_\text{Ks}-2I_\text{NaK}+I_\text{pK}+I_\mathrm{stim}}{V_cF},
\end{align*}
where $V_c$ denotes the cytoplasmic volume and $F$ the Faraday constant ($[K]_i$ is depending on the different potassium currents). For this equation it seems rather difficult or even impossible to find a steady state for which the continuation algorithm in \texttt{CL\_Matcont} converges. To get around this issue, we removed this differential equation and set $[K]_i$ constant equal to the initial concentration $[K]_i$ of the TP06 model, $[K]_i=138~mM$. 

\begin{figure}[h]
\centering
\includegraphics[width=0.8\columnwidth]{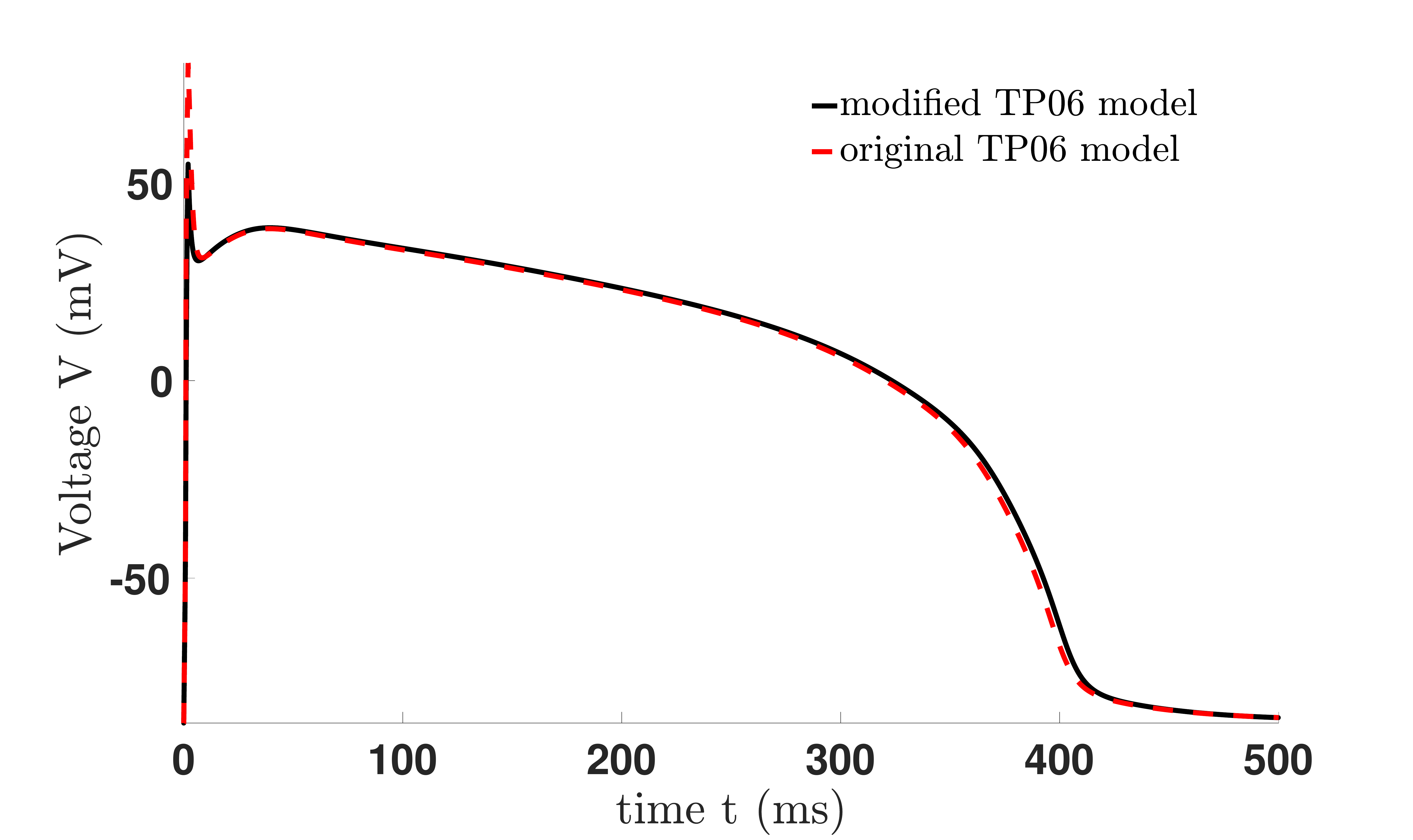}
\caption{Comparison of the trajectories of the original TP06 model with the modified one with $I_\text{stim}=52\frac{pA}{pF}$ applied for $2~ms$.}\label{fig:comparison_52}
\end{figure}

Comparing the resulting trajectories of both models, we see that they do not fit perfectly, cf. Figure~\ref{fig:comparison_52}.
To balance this, we applied an external stimulus to our modified model of $I_\text{stim}=71.5~pA/pF$ for $2~ms$ to produce Figure~\ref{fig:comparison}(a) and a slightly different one, $I_\text{stim}=73.5~pA/pF$ for $2~ms$, to gain Figure~\ref{fig:comparison}(b). Note that the trajectories are now identical. We mention that the applied external stimulus affects the initial value/starting point of the system and it has an influence on the trajectory, cf.~\cite{ES}.

The advantage of these modest alterations of an already solid model is that we are now able to investigate complex dynamics by means of bifurcation analysis. We again like to stress that, as illustrated in the figures, the reformulated model is a fairly good approximation of the original one. Thus, we can to some extent extract from the analysis fairly precise predictions of the dynamics of the TP06 model~\cite{TP06}.  

\section{Potassium dynamics of the modified model}
In the previous section, we saw that replacing the ODE for the intracellular potassium ion concentration $[K]_i$ by $[K]_i\equiv138~mM$ and the reformulation of the sodium current $I_\text{Na}$, do not affect the dynamics of the modified model too much and can be at least compensated by a stronger external stimulus. To test potential limitations of this approach we now briefly study the potassium dynamics of the modified model. 

\begin{figure}[h]
\centering
\includegraphics[width=\columnwidth]{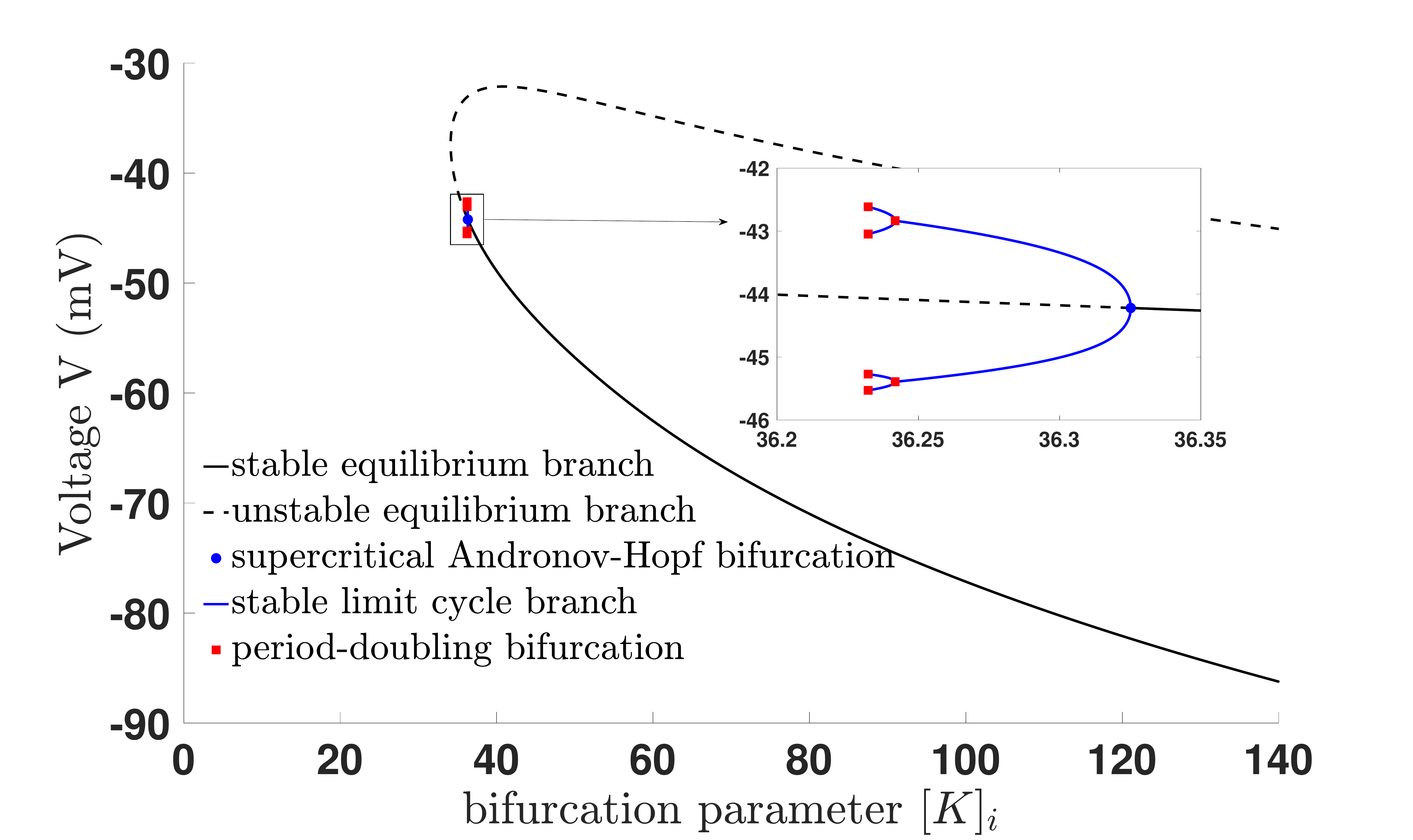}
\caption{Potassium dynamics: Bifurcation diagram for the modified TP06 model with intracellular potassium ion concentration $[K]_i$ as bifurcation parameter (projection onto the $([K]_i,V)$-plane).}\label{fig:bif_plot_K_i}
\end{figure}

To this end, we choose $[K]_i$ as bifurcation parameter and we start from $[K]_i=138~mM$ to gain the bifurcation diagram in Figure~\ref{fig:bif_plot_K_i}.
It contains a stable and unstable equilibrium branch (black solid and dashed line), respectively, and a supercritical Andronov--Hopf (AH) bifurcation (blue dot), from which a stable limit cycle branch (blue line) bifurcates. Furthermore, via a period doubling (PD) bifurcation (red square) the first limit cycle branch loses stability and a second stable limit cycle branch containing a second PD bifurcation bifurcates. 

\begin{figure*}[h]
\centering
\begin{overpic}[width=0.8\textwidth]{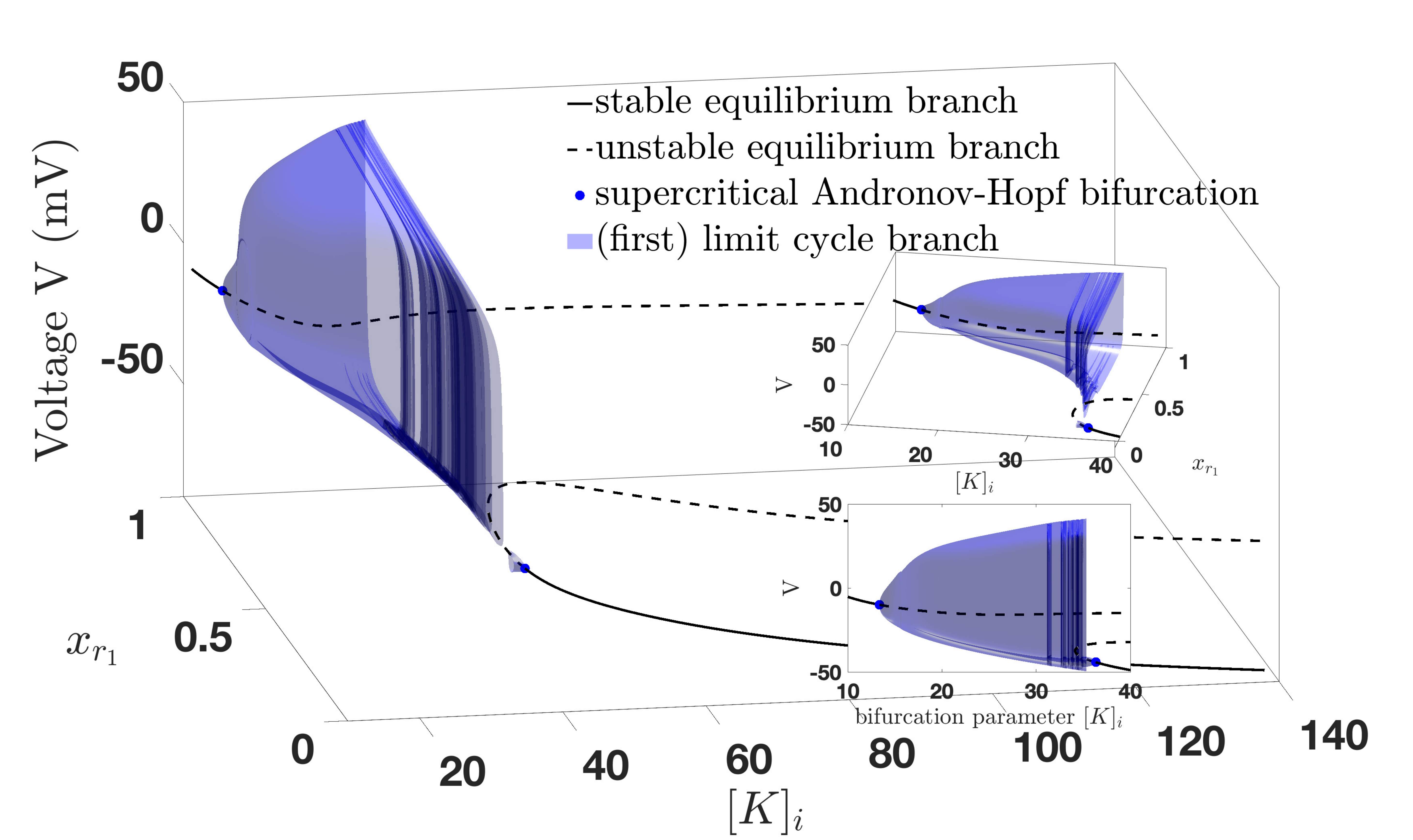}
\put(61,35){\large a)}
\put(61,22){\large b)}
\end{overpic}
\caption{Potassium dynamics: Bifurcation diagram for the modified TP06 model with $[K]_i$ as bifurcation parameter including both AH bifurcation (projection onto the $([K]_i,x_{r_1},V)$-space). a) Zoom and rotation of the bifurcation diagram. b) Zoom and projection of the bifurcation diagram onto the $([K]_i,V)$-plane.}\label{fig:bif_plot_K_i_1}
\end{figure*}

Starting from $[K]_i=10~mM$ we get a second equilibrium branch with a second supercritical AH bifurcation, cf. Figure~\ref{fig:bif_plot_K_i_1}. 
Thus, the system exhibits two supercritical AH bifurcations, one for $[K]_i\approx 13.3562~mM$ and one for $[K]_i\approx 36.3252~mM$. From both AH bifurcations a stable limit cycle branch bifurcates. Indeed, the dominating dynamics are due to the one at $[K]_i\approx 13.3562~mM$ and its limit cycle branch, respectively. 

\begin{figure*}[h]
\centering
\begin{overpic}[width=1\columnwidth]{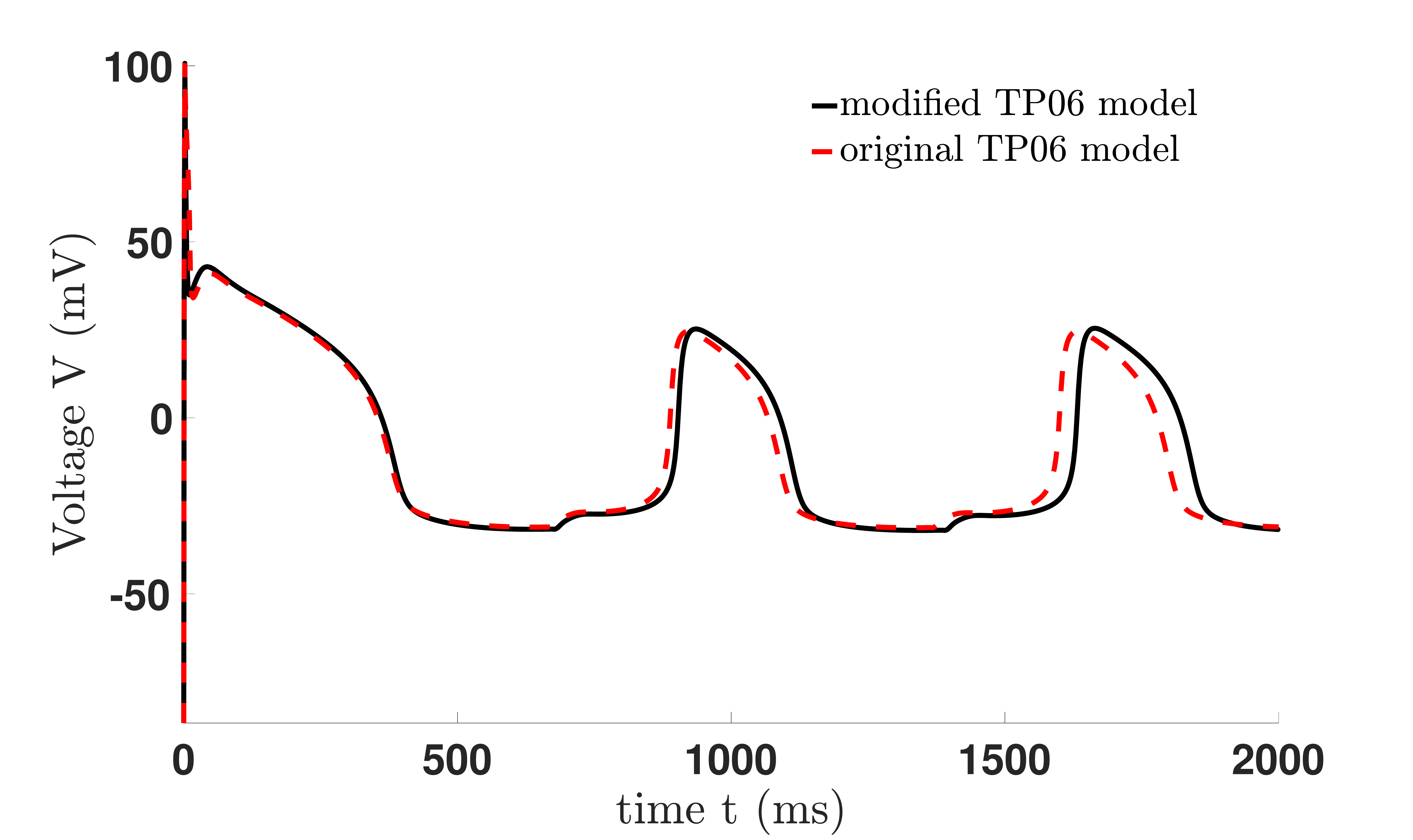}
\put(0,50){(a)}
\end{overpic}
\begin{overpic}[width=1\columnwidth]{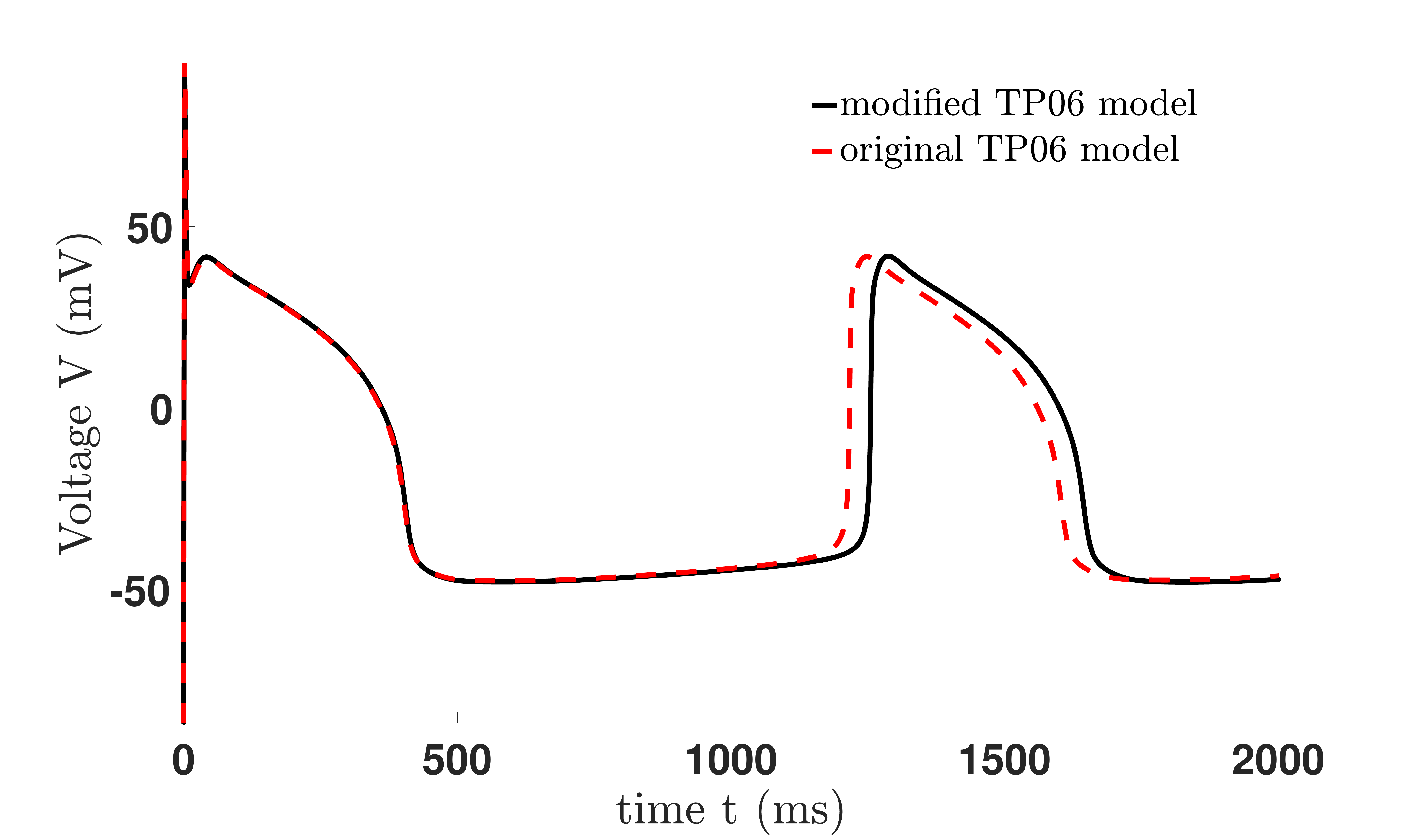}
\put(0,50){(b)}
\end{overpic}
\caption{Comparison of the trajectories of the original TP06 model (initial value $[K]_i=20~mM$ or $[K]_i=35~mM$) with the modified TP06 model with $I_\text{stim}=73\frac{pA}{pF}$ applied for $2~ms$ and $[K]_i=20~mM$ (a) or $[K]_i=35~mM$ (b).}\label{fig:comparison_K_i}
\end{figure*}

Notice that both limit cycle branches collide with the unstable equilibrium branch and terminate there. In addition, it is obvious from Figure~\ref{fig:bif_plot_K_i_1} that the modified TP06 model may exhibit periodic oscillatory pattern for a reduced intracellular potassium concentration $[K]_i$. Figure~\ref{fig:bif_plot_K_i_1} also indicates that reducing $[K]_i$ at least determines whether the AP reaches it resting potential or not. For instance, $[K]_i=40~mM$ implies that the trajectory has a stable equilibrium at approximately $-50~mV$. Therefore, the resting potential cannot be reached, cf. Figure~\ref{fig:comparison_K_i}.

Comparing these findings with the original TP06 model by choosing the $[K]_i$ value of the modified TP06 model as the initial intracellular potassium concentration of the original TP06 model, we see that our bifurcation analysis yields also a fairly good predication for the dynamics of the original TP06 model, cf. Figure~\ref{fig:comparison_K_i}. More precisely, a reduction of the initial intracellular potassium concentration for the original TP06 model may lead to periodic oscillatory pattern. 

This shows that fixing and varying $[K]_i$ leads to similar dynamics as changing the initial intracellular potassium concentration of the 
TP06 model. In both cases, the trajectories do not reach the resting potential and no further initial stimulus is needed to produce a change in the membrane potential, since periodic patterns occur, some sort of ventricular tachyarrhythmias~\cite{WEISS20101891,Hypokalemia,Hypokalemia1}. 
Due to the fact that the resting phase is designated by high potassium 
currents, it is not surprising that the reduction of $[K]_i$ may inhibit the AP to reach its resting \mbox{potential.}
\indent Our analysis shows that considering a fixed intracellular potassium concentration is not a significant limitation and yields comparable results as in the TP06 model choosing the fixed $[K]_i$ value as initial value. The modification of the sodium current $I_\mathrm{Na}$ also modestly modify the behaviour of the system. Therefore, these modifications are reasonable and one can keep $[K]_i\equiv138~mM$ fixed.

\section{Reduction of the rapid and slow potassium currents}
\begin{figure}[h]
\subfigure[$G_\mathrm{K_s}= 0.02505\frac{nS}{pF}$.]{\includegraphics[width=0.5\columnwidth]{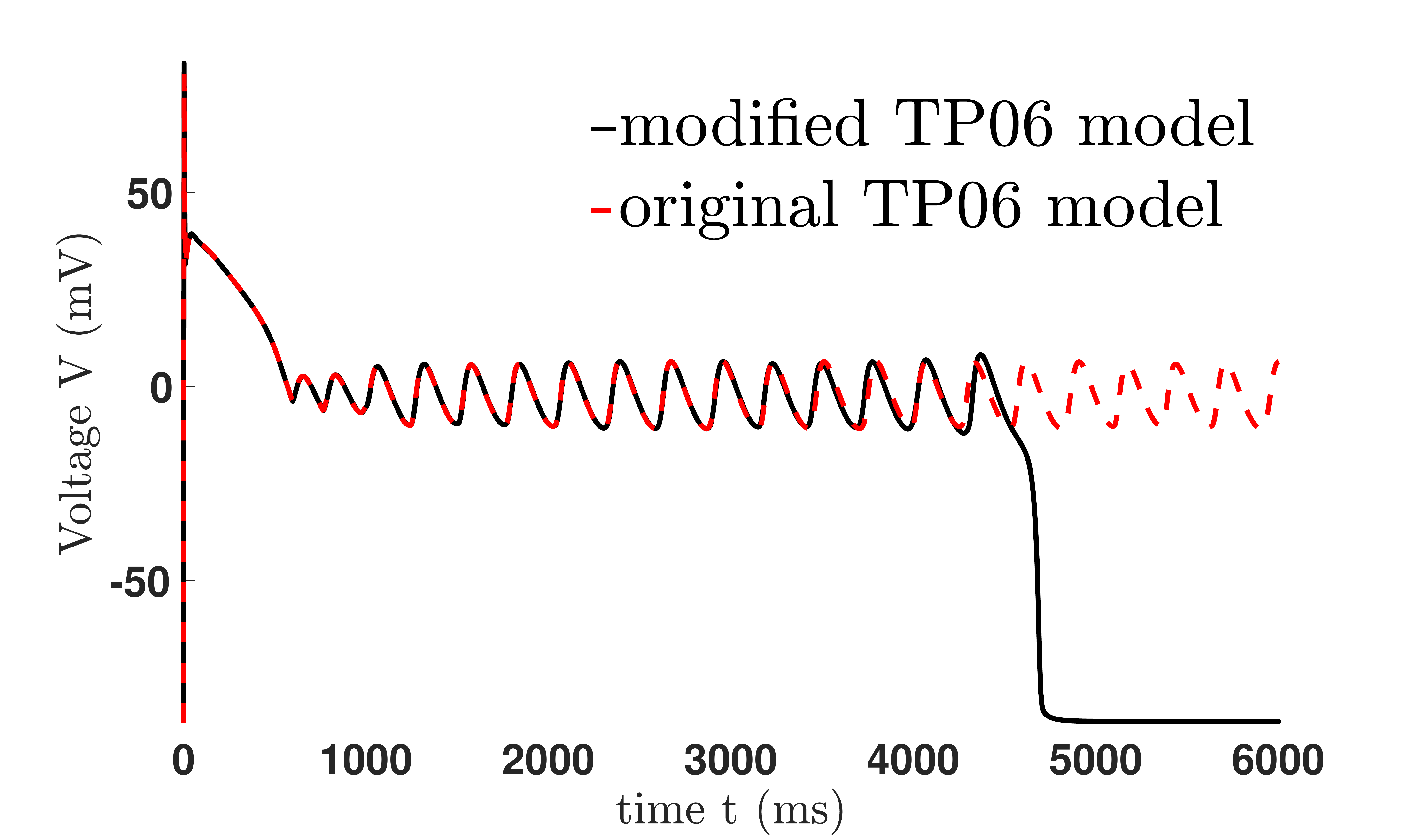}}\subfigure[$G_\mathrm{K_s}= 0.02505\frac{nS}{pF}$.]{\includegraphics[width=0.5\columnwidth]{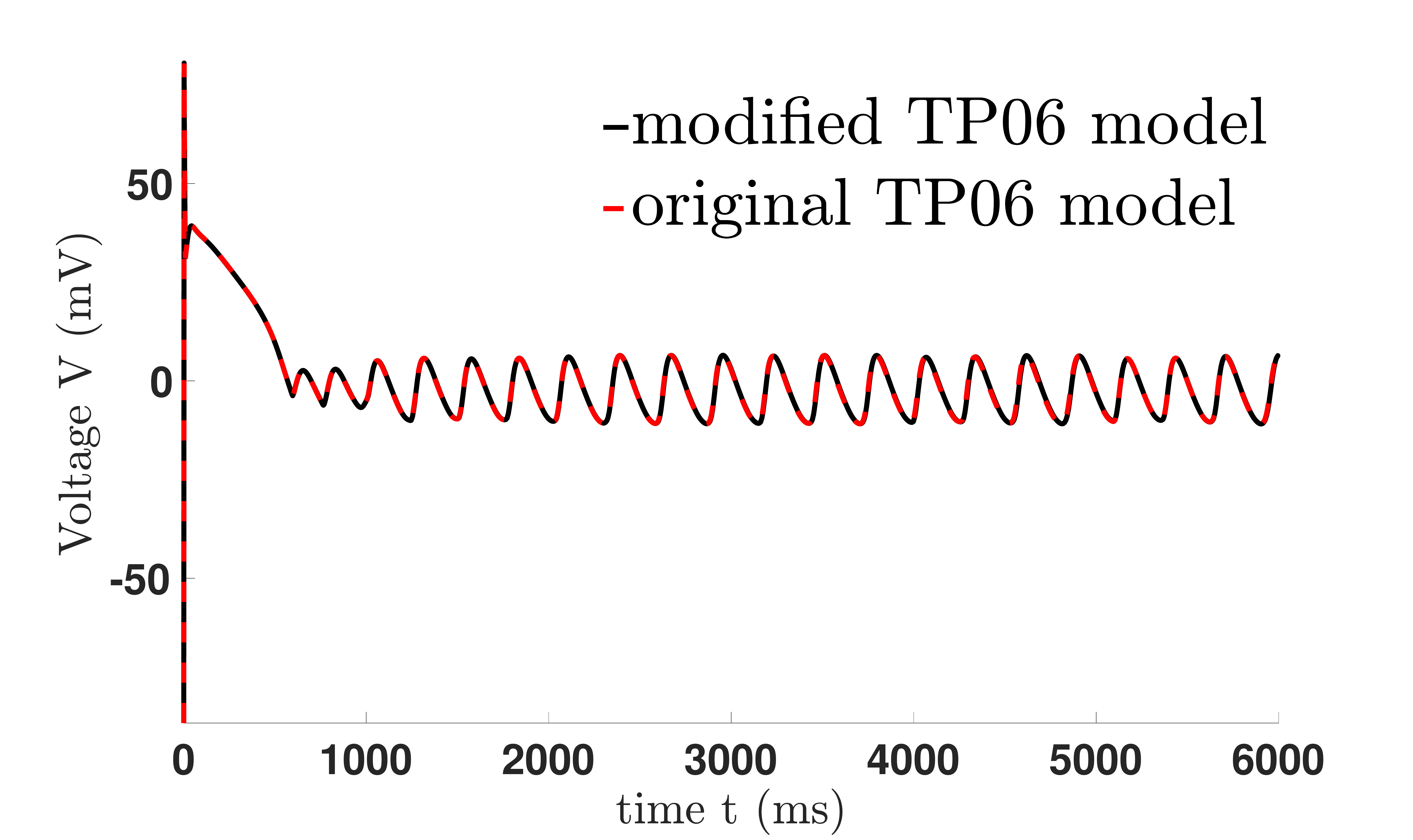}}
\\
\subfigure[$G_\mathrm{K_s}= 0.028\frac{nS}{pF}$.]{\includegraphics[width=0.5\columnwidth]{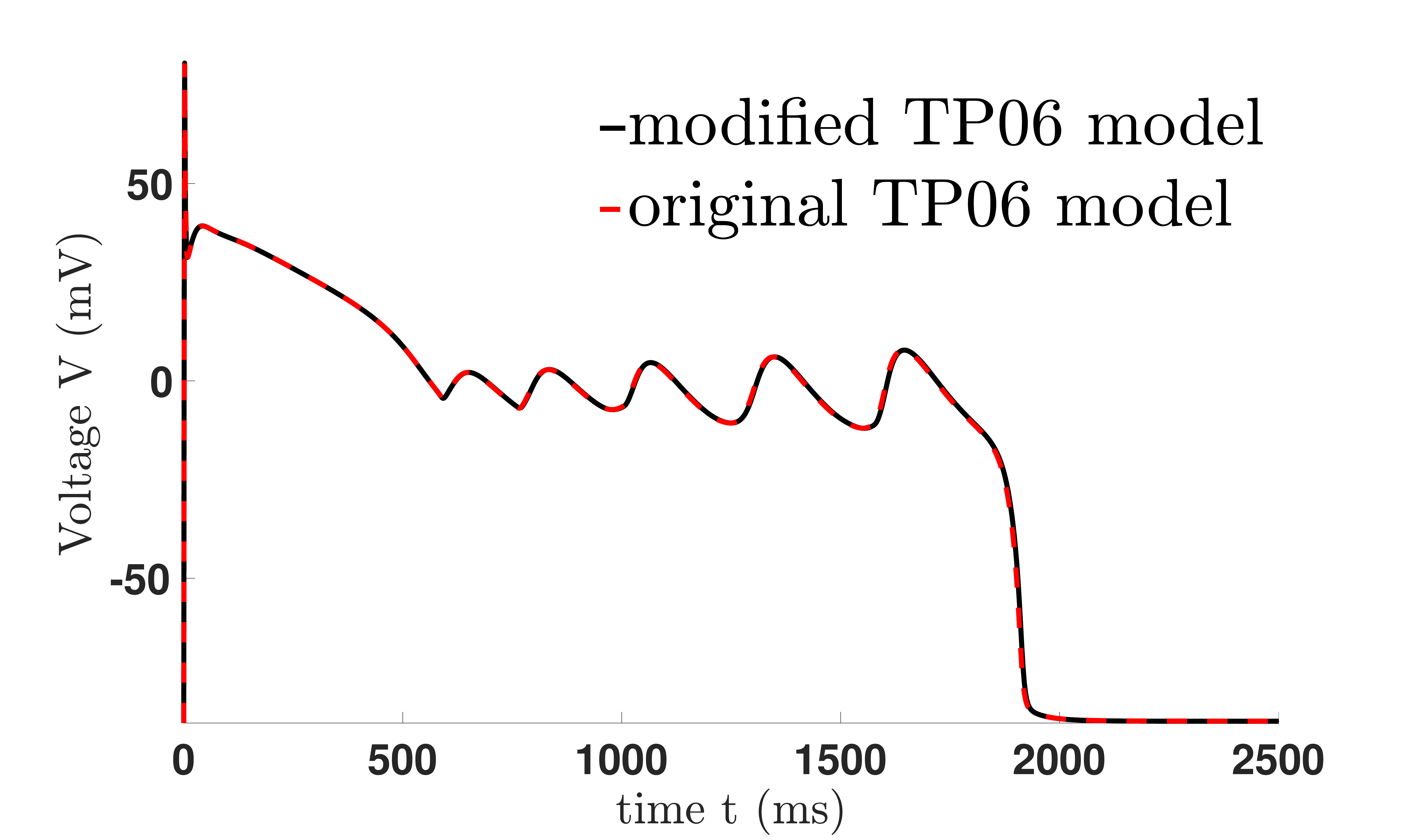}}\subfigure[$G_\mathrm{K_s}= 0.04\frac{nS}{pF}$.]{\includegraphics[width=0.5\columnwidth]{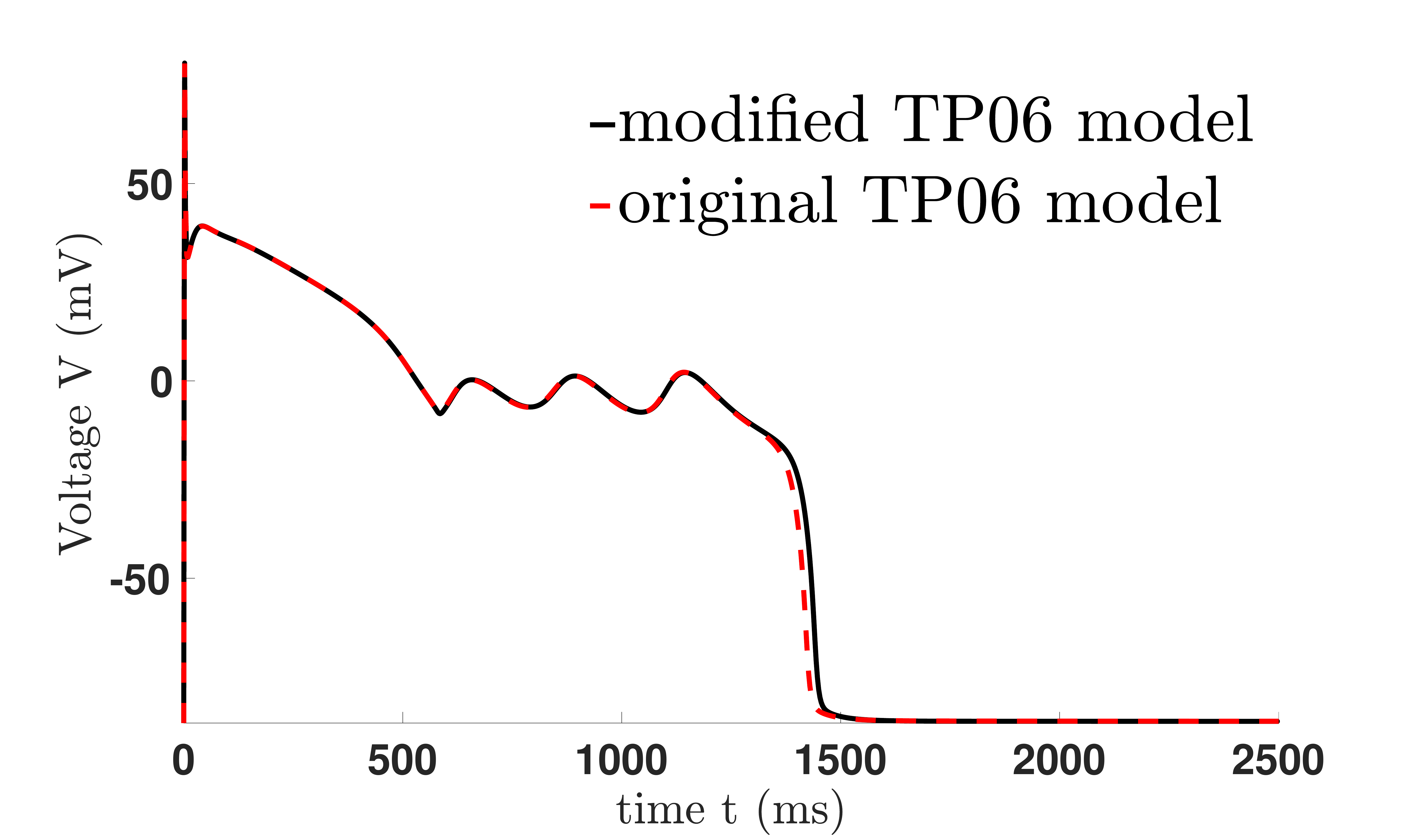}}
\caption{Comparison of the trajectories of the modified TP06 and the TP06 model with different reduced $G_\mathrm{K_s}$ value, where $G_\mathrm{K_r}= 0.153~nS/pF$ is the value from~\citep{TP06}. Furthermore, in (b) - (d) a stimulus  $I_\text{stim}=71.5~pA/pF7$ over $2~ms$ was applied, while in (a) $I_\text{stim}=73.5~pA/pF$ over $2~ms$.}\label{fig:comparison_reduced_GKs}
\end{figure}

Due to the common knowledge that a reduction of the potassium currents, mainly the rapid one $I_\mathrm{K_r}$, may lead to certain cardiac arrhythmia like EADs, it is natural to investigate the dynamics of the modified TP06 model by means of bifurcation analysis choosing one or more of the current specific conductances $G_\mathrm{K_r}$ and/or $G_\mathrm{K_s}$ as bifurcation parameter(s). Starting with the dynamics of the modified TP06 model with respect to $G_\mathrm{K_r}$ it turns out that the system has only for negative values of $G_\mathrm{K_r}$ an AH bifurcation and the limit cycle branch remains also on the negative part of the $G_\mathrm{K_r}$-axis. Furthermore, no EADs (via a reduced rapid potassium current) appear, neither for the modified TP06 model nor for the original one, which corresponds also to the findings in \cite{SpatialPatterns}. However, EADs may appear via a reduced slow potassium current $I_\mathrm{K_s}$ and as a combination of a reduction of the fast and slow potassium currents. The potential effect of a reduced slow potassium current $I_\mathrm{K_s}$ is illustrated in Figure \ref{fig:comparison_reduced_GKs}.

\begin{figure}[h]
\centering
\subfigure[Bifurcation diagram 2D.]{\includegraphics[width=\columnwidth]{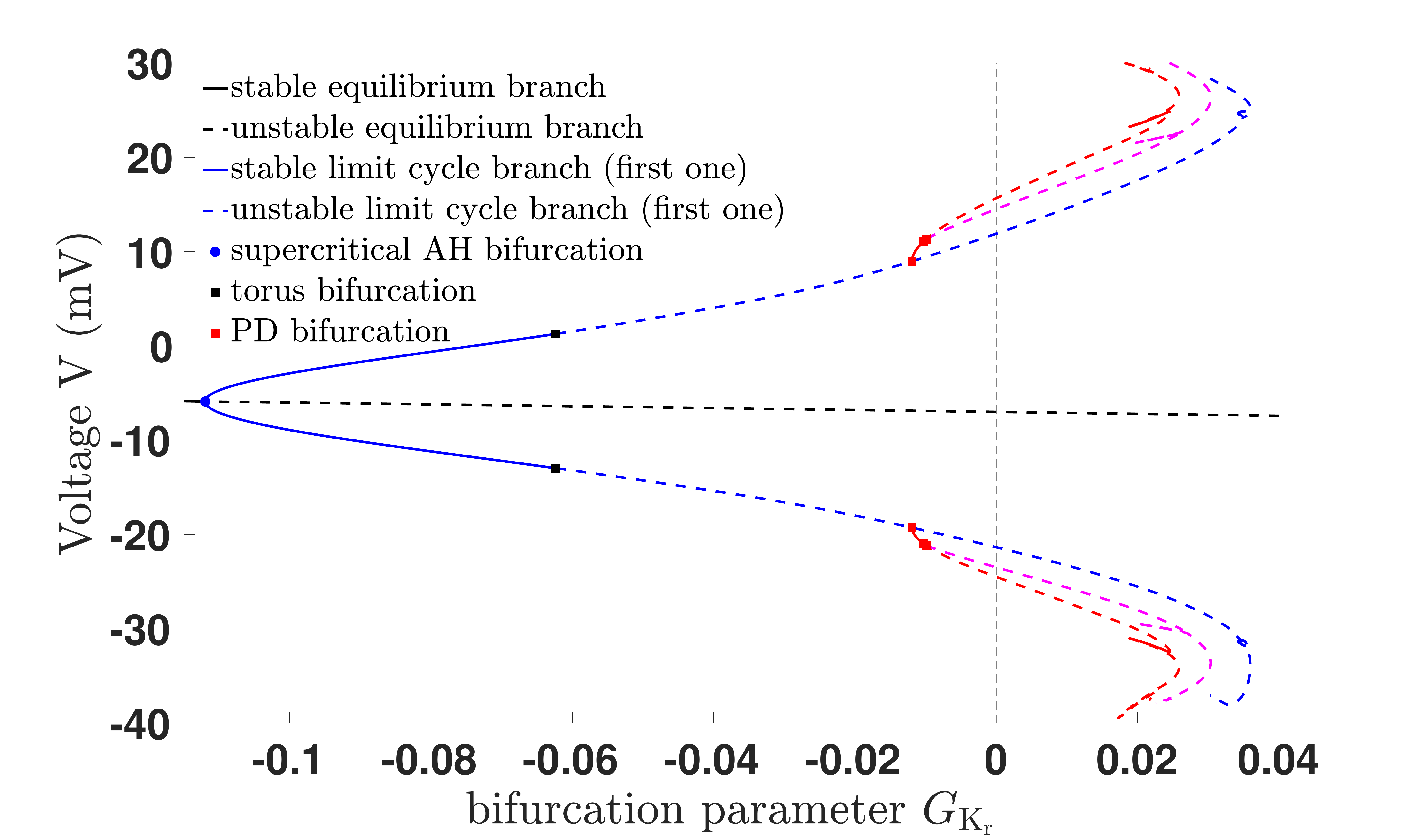}}
\subfigure[Bifurcation diagram 3D.]{\includegraphics[width=\columnwidth]{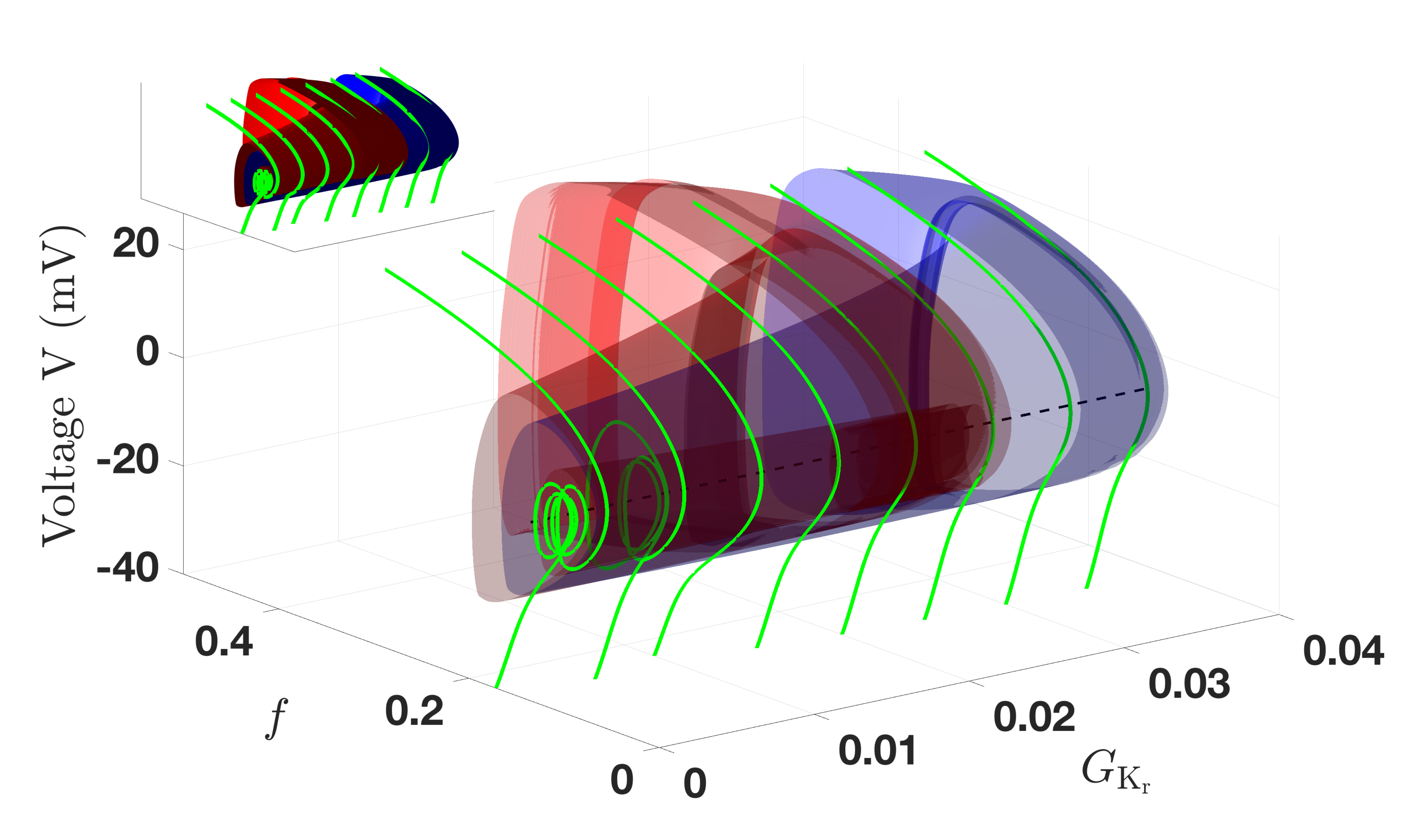}}
\caption{Bifurcation diagram of the modified model of~\citep{TP06} with a $75\%$ blockade of the slow potassium current $I_\mathrm{K_s}$, using $G_\mathrm{K_r}$ as bifurcation parameter: (a) Projection onto the $(G_\mathrm{K_r},V)$-plane. (b) Projection onto the $(G_\mathrm{K_r},f,V)$-space.
}\label{fig:bif_combination}
\end{figure}

Therefore, we focus our study on the dynamics of the modified TP06 model with respect to a combination of $G_\mathrm{K_s}$ and $G_\mathrm{K_r}$. As one example for the occurrence of EADs via a combination of reduced $G_\mathrm{K_s}$ and $G_\mathrm{K_r}$, we consider a $75\%$ blockade of the slow potassium current $I_\mathrm{K_s}$, i.e. $\bar{G}_\mathrm{K_s}=0.25\cdot G_\mathrm{K_s}$, and we choose the conductance $G_\mathrm{K_r}$ as bifurcation parameter, cf. Figure~\ref{fig:bif_combination}(a) \& (b). Remember that our model is a 17 dimensional ODE system, thus we have a 17 dimensional phase space plus an additional dimension due to the bifurcation parameter. Therefore, the visualisation is tricky and we only present certain projections, which illustrate the dynamics well. 
In Figure~\ref{fig:bif_combination}(a) the equilibrium curve changes stability via a supercritical AH bifurcation (at $G_\mathrm{K_r}\approx -0.1120\frac{nS}{pF}$) and a limit cycle branch bifurcates, which changes stability at a torus bifurcation (at $G_\mathrm{K_r}\approx -0.0623\frac{nS}{pF}$), cf. \citep{EMS}, and exhibits a PD bifurcation (at $G_\mathrm{K_r}\approx -0.0119\frac{nS}{pF}$). Notice that the AH bifurcation lies outside a physiological relevant range as well as the torus bifurcation and the PD bifurcations, cf. Figure~\ref{fig:bif_combination}(a). However, parts of the limit cycle branches are in a plausible range and determine the dynamics of the system. From the first PD bifurcation a further limit cycle branch bifurcates and contains a second PD bifurcation (at $G_\mathrm{K_r}\approx -0.0103\frac{nS}{pF}$), where a third limit cycle branch bifurcates with a third PD bifurcation (at $G_\mathrm{K_r}\approx -0.0099\frac{nS}{pF}$), cf. Figure~\ref{fig:bif_combination}(a). Such a PD cascade is an indication of the occurrence of EADs, cf.~\citep{ES,AE_control,AE_MMOs}. Although the analysis indicates that EADs appear, whether they actually do is also determined by the external stimulus.

This is visible in Figure~\ref{fig:bif_combination}(b), where we include only the first two limit cycle branches and consider only a physiological relevant range of $G_{\text{K}_r}$. Moreover, we include several trajectories of the modified system for different values of $G_{\text{K}_r}$ (green lines). Although, the autonomous system exhibits a PD cascade, EADs only appear if the trajectory of the non-autonomous ODE system can enter the inside of the (unstable) limit cycle branches and can surround the unstable equilibrium curve before leaving again the basin of attraction and reaching the resting potential. This indicates how initial values and/or the initial stimulus also decide whether a trajectory reaches a certain basin of attraction and certain dynamics appear. Finally, we want to highlight that the first trajectory in Figure~\ref{fig:bif_combination}(b), the one for $G_\mathrm{K_r}=0\frac{nS}{pF}$, is the same trajectory as in Figure~\ref{fig:comparison}(b).

\begin{figure}[h]
\centering
\subfigure[Bifurcation diagram 2D.]{\includegraphics[width=\columnwidth]{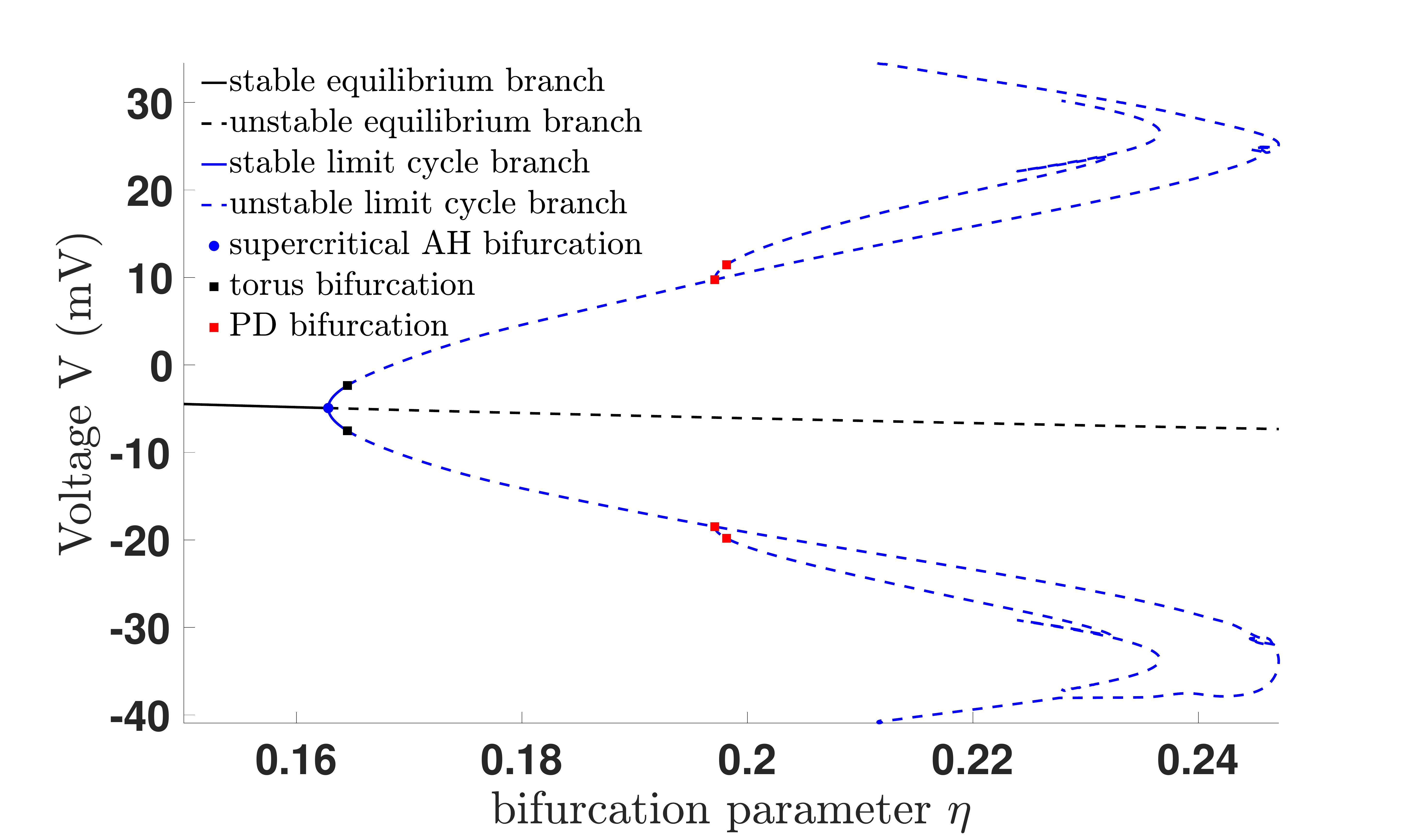}}
\subfigure[Bifurcation diagram 3D.]{\includegraphics[width=\columnwidth]{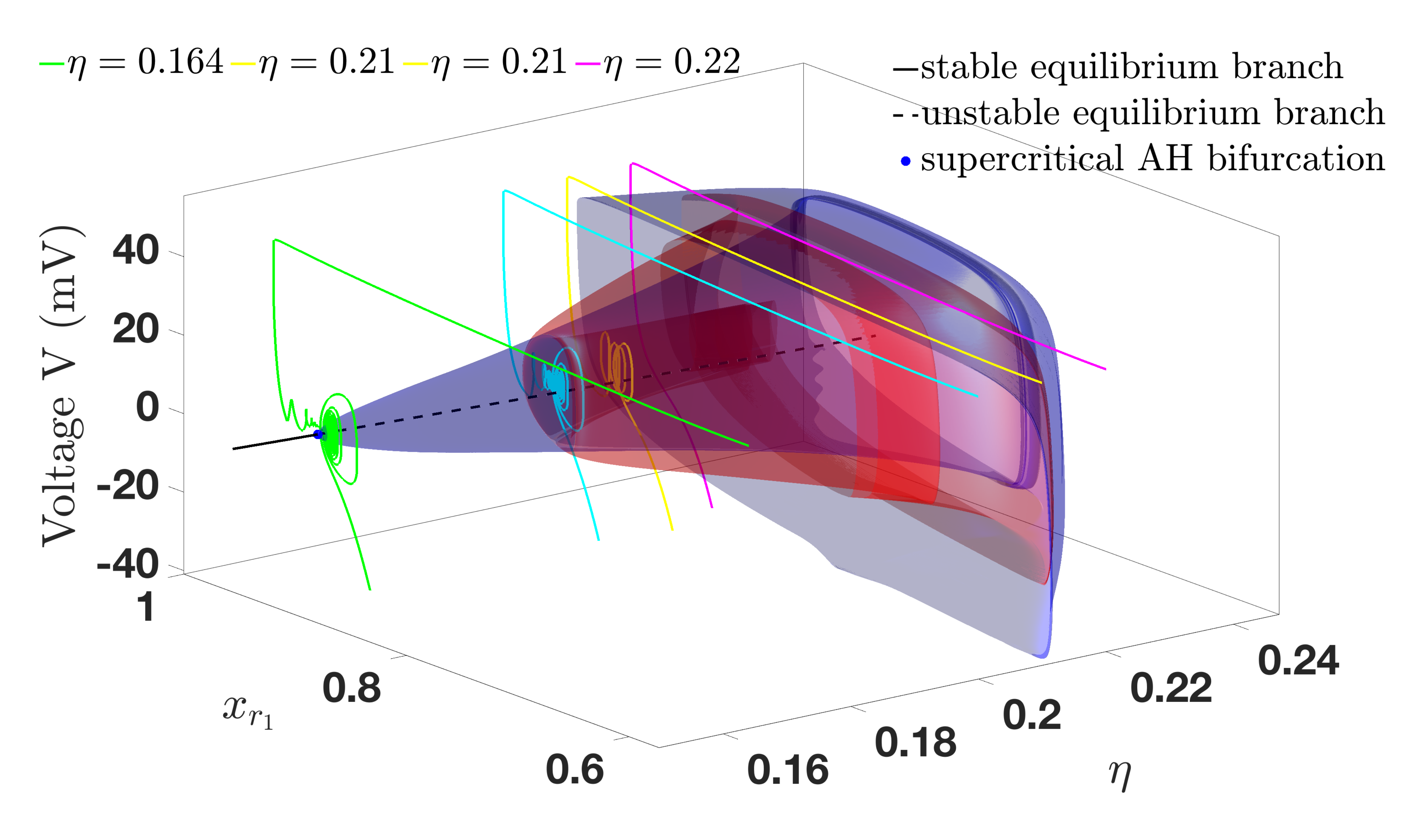}}
\caption{Bifurcation diagram of the modified model of~\citep{TP06} utilising $\eta$, cf. equation~\eqref{eta}, as bifurcation parameter: (a) Projection onto the $(\eta,V)$-plane. (b) Projection onto the $(\eta,x_{r_1},V)$-space.
}\label{fig:bif_Ks_Kr}
\end{figure}

\begin{figure}[h]
\centering
\includegraphics[width=\columnwidth]{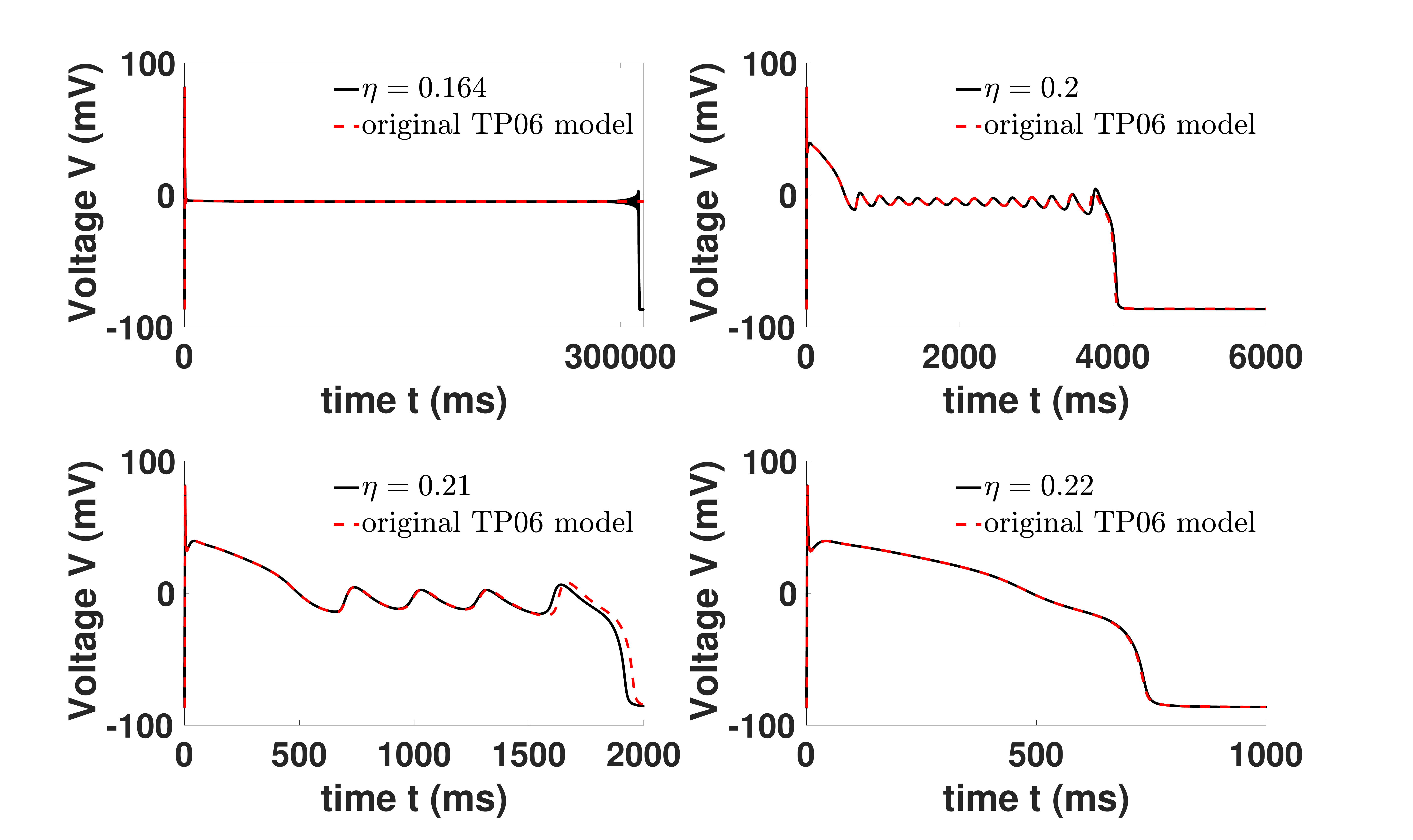}
\caption{Trajectories of the modified TP06 model for different values of $\eta$, which are included in the bifurcation diagram \ref{fig:bif_Ks_Kr}(b).}\label{fig:bif_Ks_Kr_sim}
\end{figure}

In addition, we investigate the effect of a blockade of both potassium currents by introducing a new parameter $\eta\in[0,1]$:
\begin{align}
    \bar{G}_\mathrm{K_r}=\eta \cdot G_\mathrm{K_r}\quad\text{and}\quad\bar{G}_\mathrm{K_s}=\eta \cdot G_\mathrm{K_s}.\label{eta}
\end{align}
This enables the study how a potential blockade of the fast and slow potassium current at the same time affects the dynamics of the modified TP06 model. This yields the bifurcation diagram in Figure~\ref{fig:bif_Ks_Kr} using $\eta$ as bifurcation parameter reflecting the percentage of the blockade. Similar to Figure~\ref{fig:bif_combination} the modified TP06 model exhibits a supercritical AH bifurcation for $\eta \approx0.1628$, where a stable limit cycle branch bifurcates, which loses stability at a torus bifurcation ($\eta \approx 0.1645$) and from the first PD bifurcation ($\eta\approx 0.1971$) a second unstable limit cycle branch bifurcates. In this situation the dynamics of our interest, i.e. EADs, can be connected with the torus bifurcation (EADs with an extreme prolongation of the duration, but only for the modified TP06 model) and with the PD bifurcations (EADs with a ``normal'' prolongation of the duration for both models), cf. Figure~\ref{fig:bif_Ks_Kr}(b) and Figure~\ref{fig:bif_Ks_Kr_sim}. We see that normal AP appears already as soon as the first limit cycle branch (with an unstable, repelling branch) covers the complex dynamics of the autonomous system. We also have in a physiologically relevant region a transition between an excitable and oscillatory media close to the supercritical AH bifurcation as in~\citep{ES}, cf.~\cite{sandstede2021spiral,Dodson}.


\section{Predictions of the modified model on the dynamics of the TP06 monodomain model}

One major aim in cardiac research is to understand wave patterns like travelling waves, spiral waves, wave breakup and instabilities that might be linked to certain cardiac arrhythmia, see~\cite{sandstede2021spiral,Alonso_2016,Cherry2017,TVEITO2008141}. In this section, we show by example that complex dynamics linked to EADs found through the bifurcation analysis of the modified model, yields complex dynamics on the tissue level of the TP06 model. 

To study the behaviour of the model on macro-scale, we use the following 2D monodomain model:
\begin{align}
C_m\frac{\partial V}{\partial t} = - I_\text{ion}+I_\mathrm{stim} + \nabla\cdot\left(\begin{pmatrix}D&0\\0&D\end{pmatrix}\nabla V\right), \label{monodomain}
\end{align} 
where $C_m=1~\mu F/cm^2$ and $D=0.00154~cm^2/ms$, equipped with Neumann boundary condition. The diffusion matrix is chosen in accordance with~\cite{TP06}. Note that the mondomain model is directly linked to the bidomain model, which takes into account the anisotropy of the intra- and extracellular spaces. For a detailed description we refer to~\cite{Nielsen,computing_heart,quarteroni,PANFILOV20191,niederer2019computational,booktissue}.

The simulations of the tissue level model \eqref{monodomain} are run on a $1000\times 1000$ grid with spatial discretization $dx = 0.025$ cm with a five-point stencil for the Laplacian (as in \cite{Vandersickel1}) and with time step $dt = 0.0812$. The Rush and Larsen scheme \cite{rushlarsen} is used to integrate the gating variables in time, and Eulers method is used for the rest. 

\begin{figure}[h]
\subfigure[]{\includegraphics[width=0.49\columnwidth]{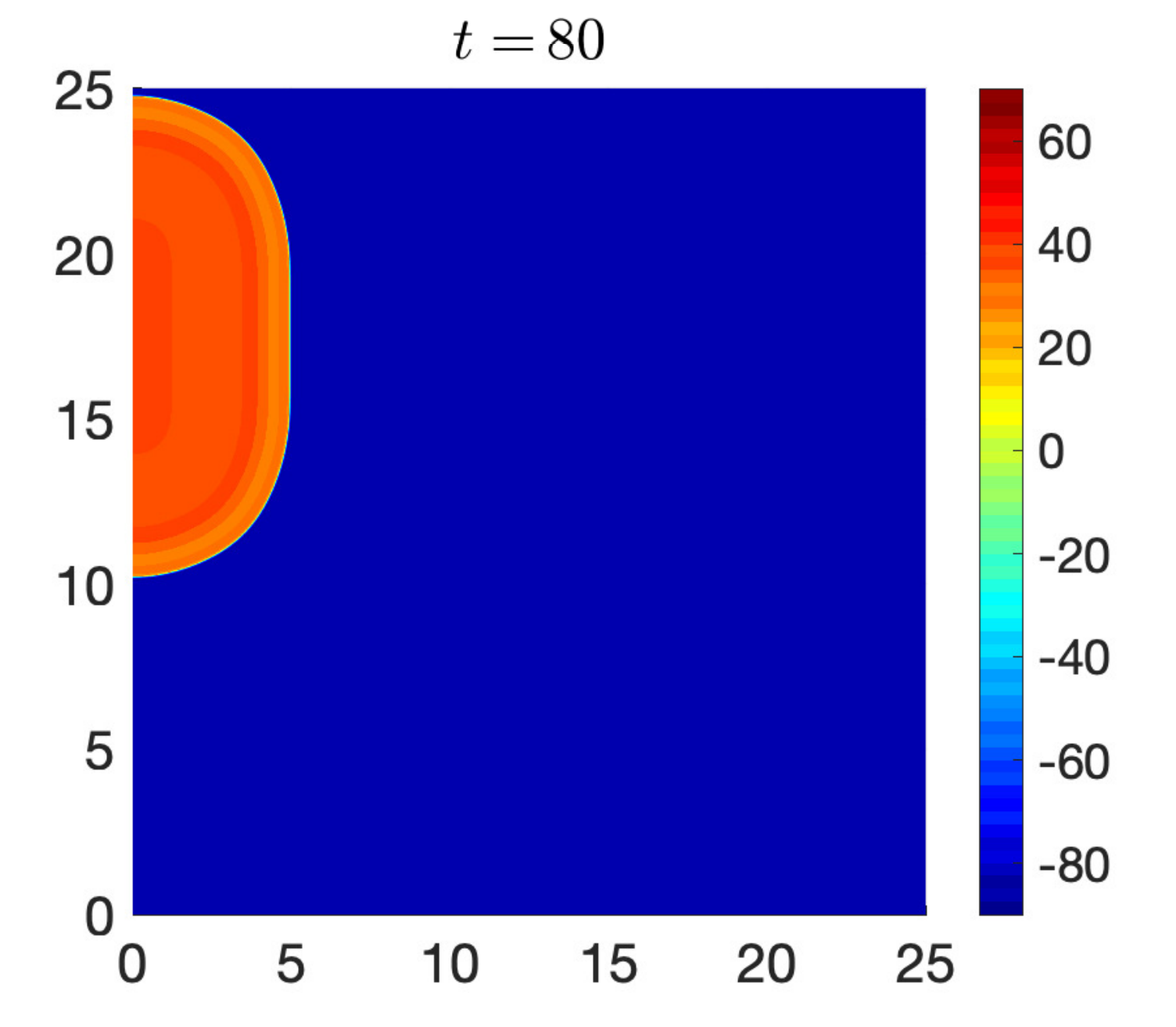}}
\subfigure[]{\includegraphics[width=0.49\columnwidth]{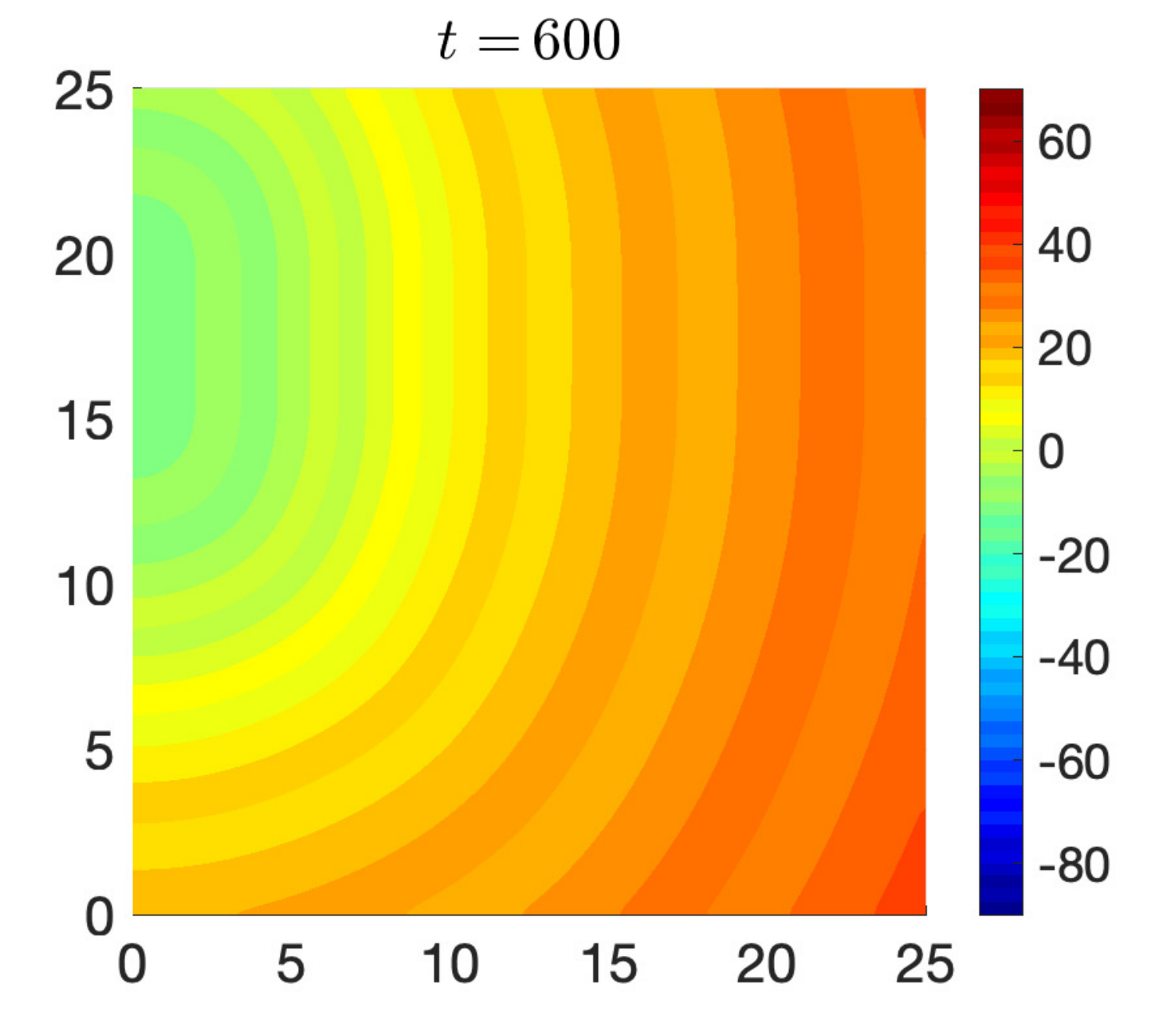}}
\\
\subfigure[]{\includegraphics[width=0.49\columnwidth]{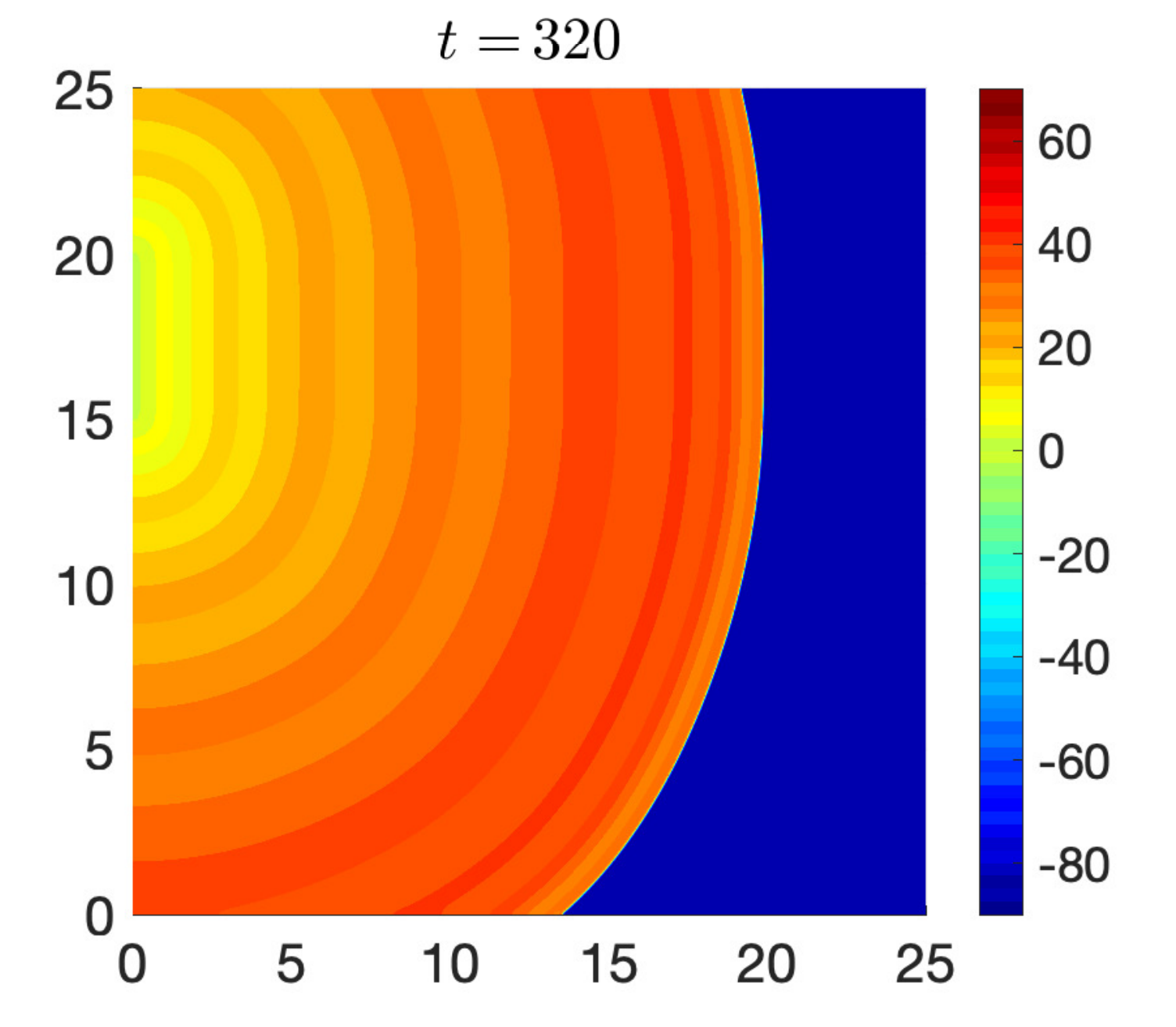}}
\subfigure[]{\includegraphics[width=0.49\columnwidth]{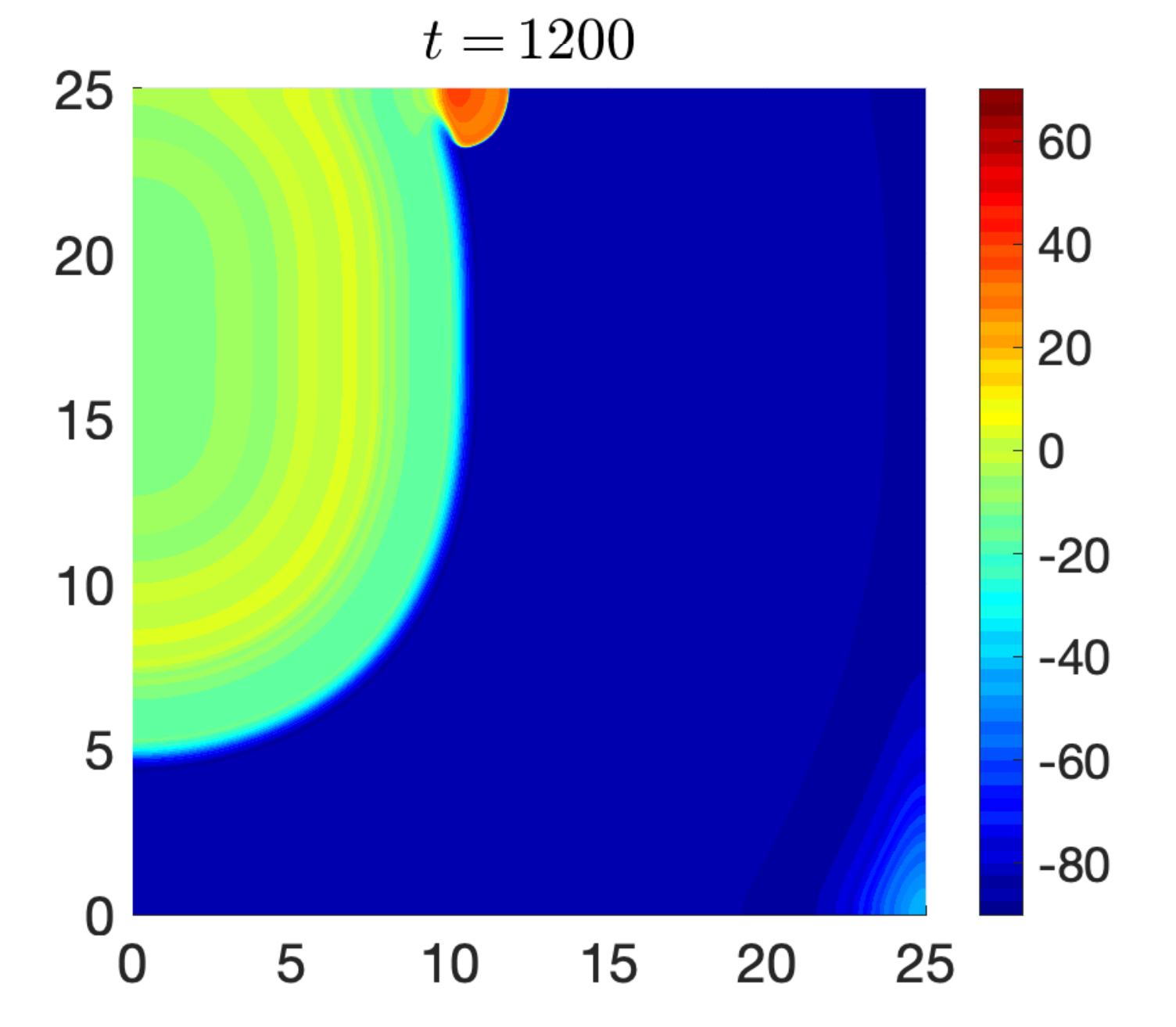}} 
\\
\subfigure[]{\includegraphics[width=0.49\columnwidth]{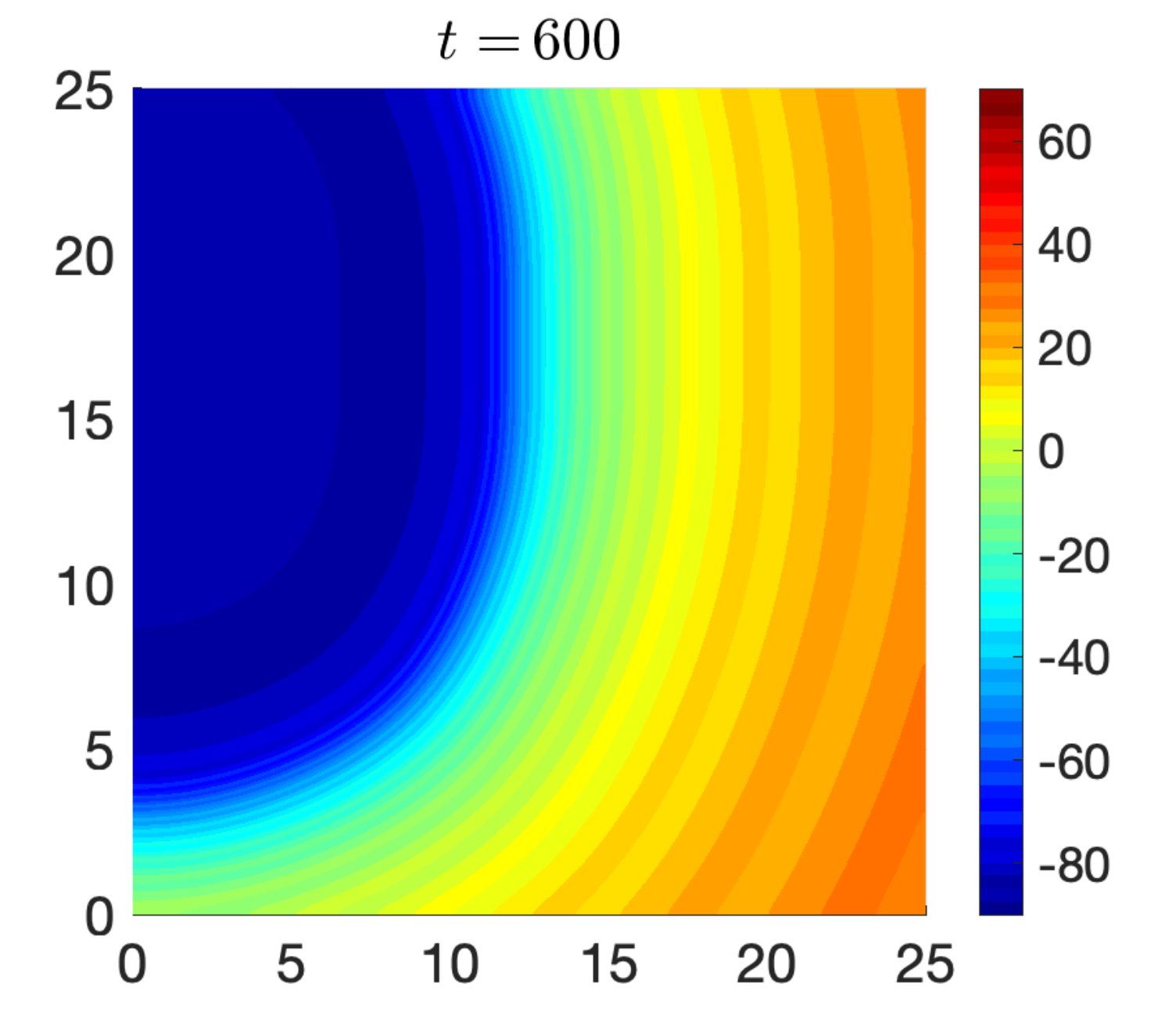}}
\subfigure[]{\includegraphics[width=0.49\columnwidth]{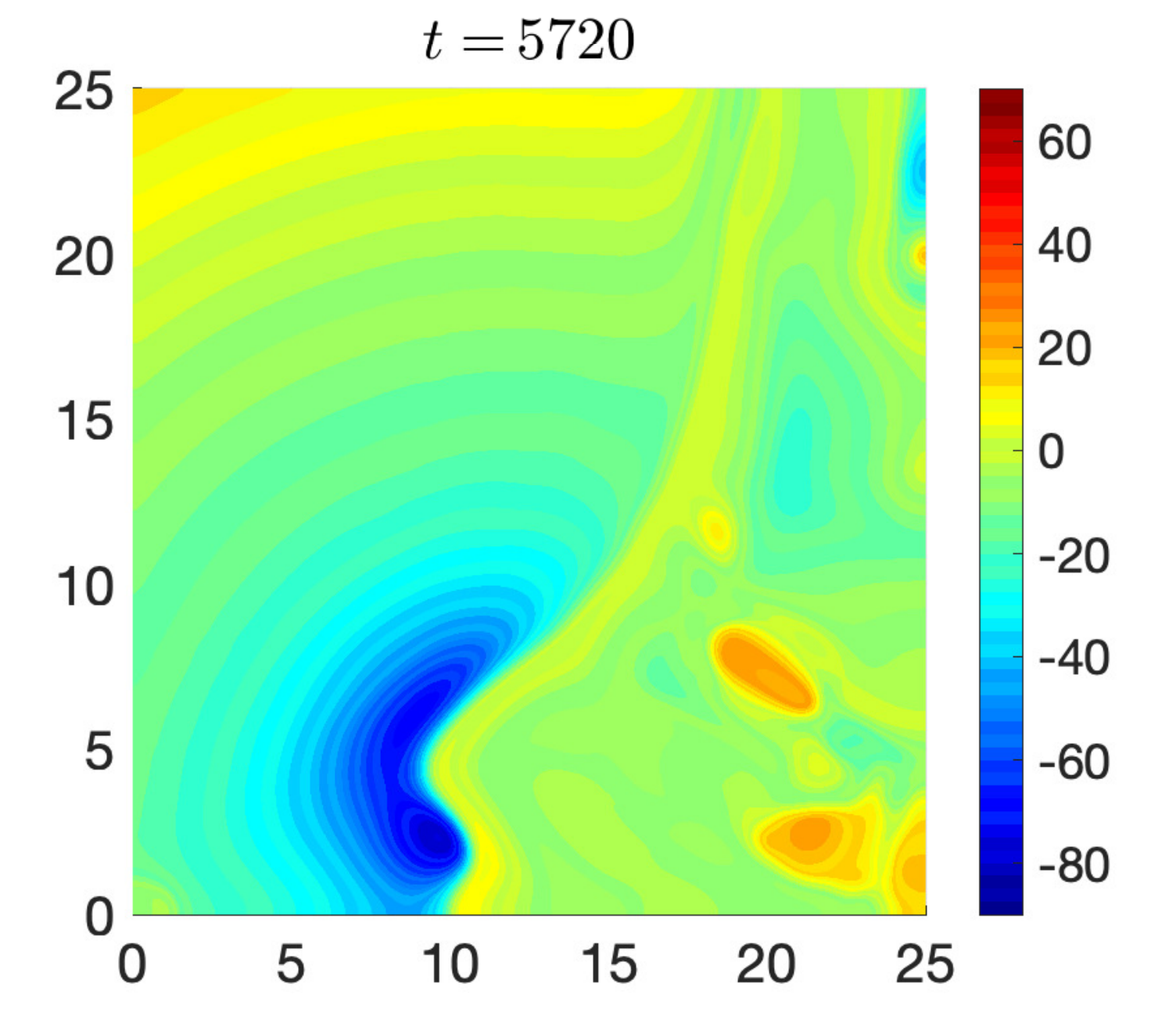}}
\caption{The monodomain model \eqref{monodomain} is simulated by applying a $6 \times 200$ cells on the left side of the domain with the two parameter settings used to produce Figure \ref{fig:comparison}. The time $t$ is measured in milliseconds. The parameters of Figure \ref{fig:comparison}(a) are used to produce the left column, ant the parameters of Figure \ref{fig:comparison}(b) are used to the right.}\label{fig:blob}
\end{figure}

In Figure \ref{fig:blob} we apply a stimulus of $I_\mathrm{stim} = 26.2 ~pA/pF$ to $6 \times 200$ cells on the left side of the domain for both the normal parameter settings of the TP06 model and for the parameter settings in Figure \ref{fig:comparison}(b).

The model \eqref{monodomain} with the normal setting, presented in the left column, produces a wave that simply crosses the computational domain. In the beginning, system \eqref{monodomain} with the other setting produces an identical wave to that of the one in the left column. Before $t=600$ms, however, this is no longer the case. After some time, the voltages at the top of the computational domain crosses the threshold (Figure \ref{fig:blob}(d)), which later results in similar voltage patterns to that of the fibrillation plots presented in \cite{Vandersickel1}. Given that the parameter setting is within the range of EADs occurring for the modified model, this is to be expected for the original TP06 model.

\begin{figure}[h]
\subfigure[]{\includegraphics[width=0.49\columnwidth]{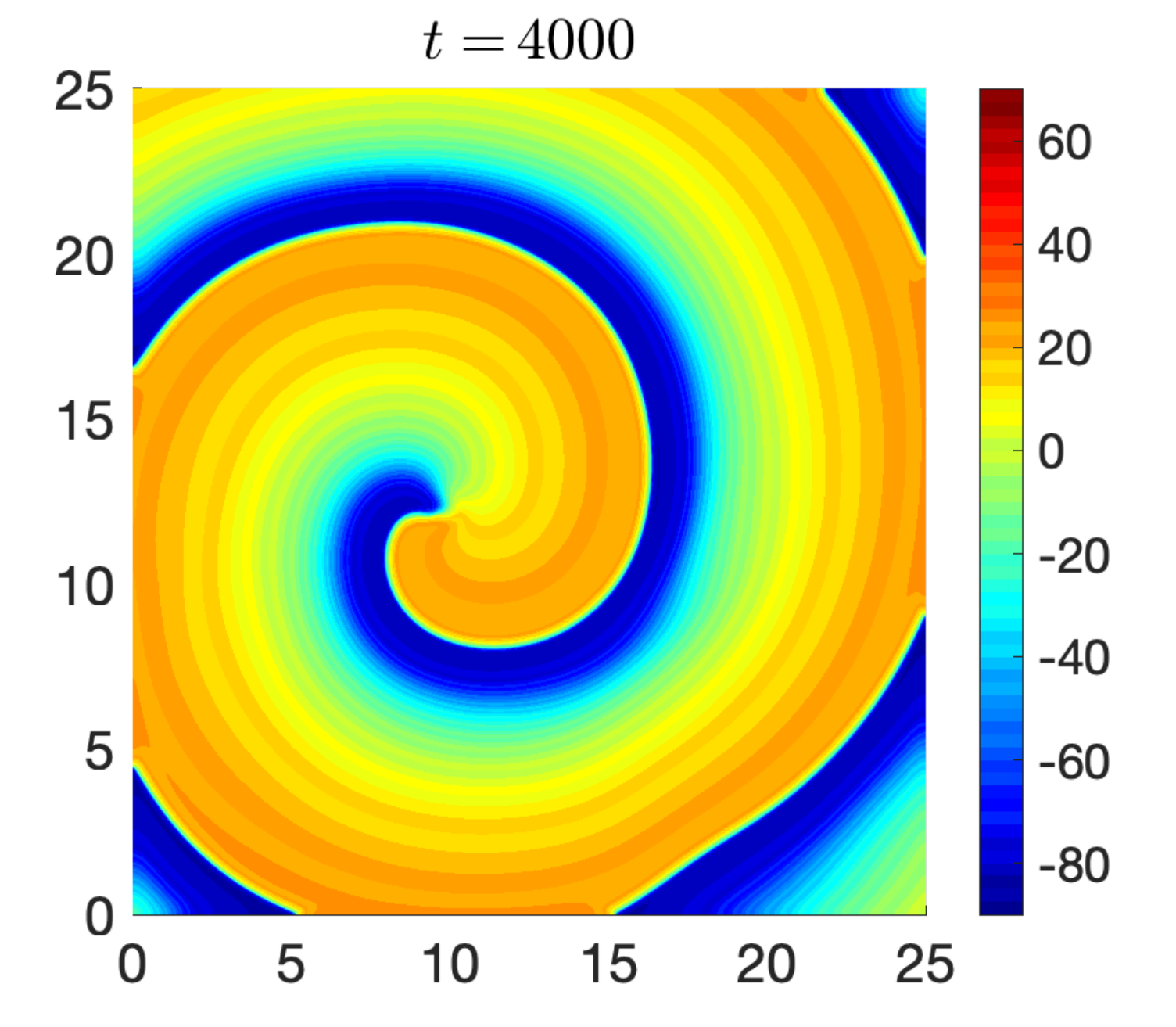}}
\subfigure[]{\includegraphics[width=0.49\columnwidth]{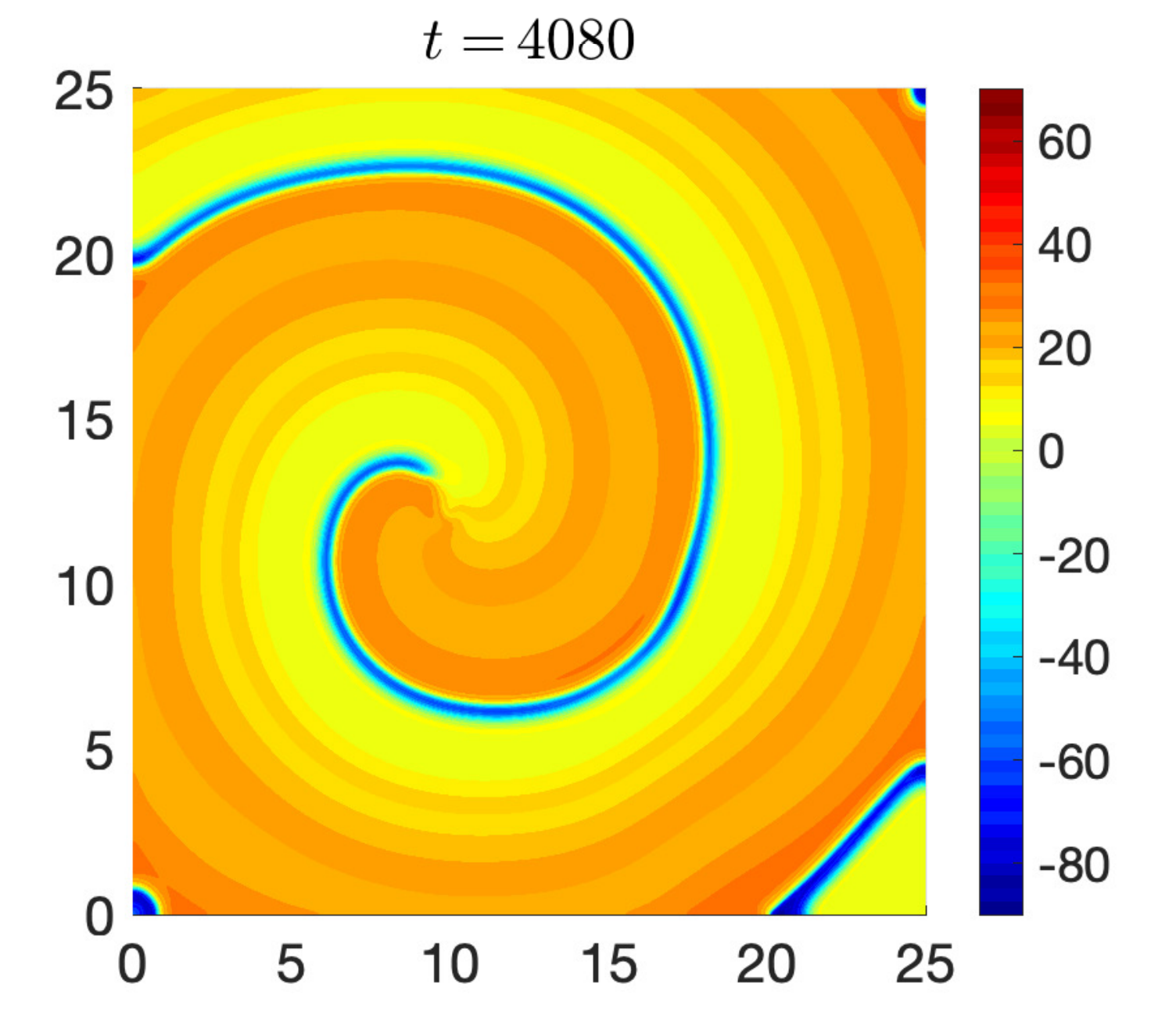}}
\\
\subfigure[]{\includegraphics[width=0.49\columnwidth]{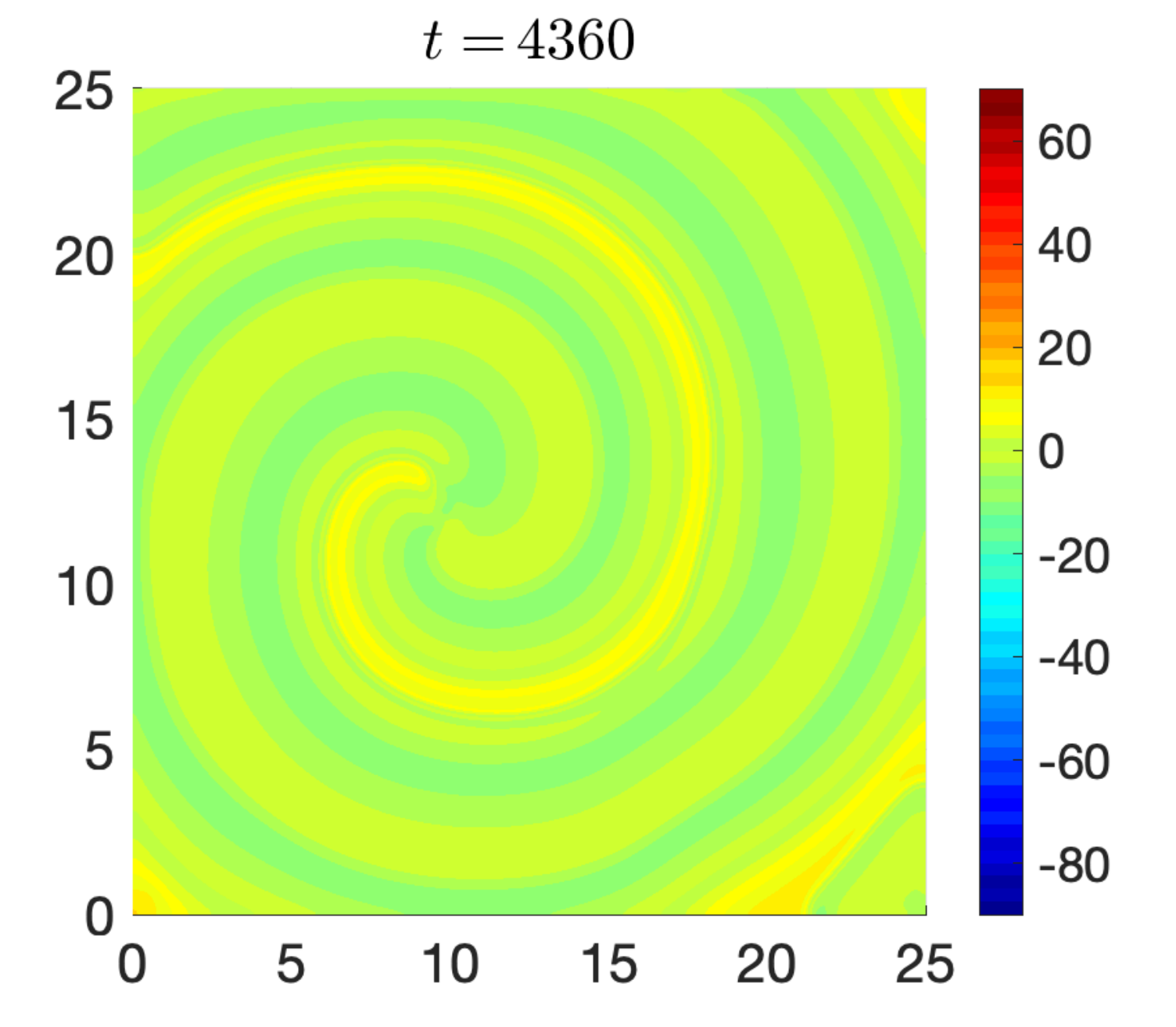}}
\subfigure[]{\includegraphics[width=0.49\columnwidth]{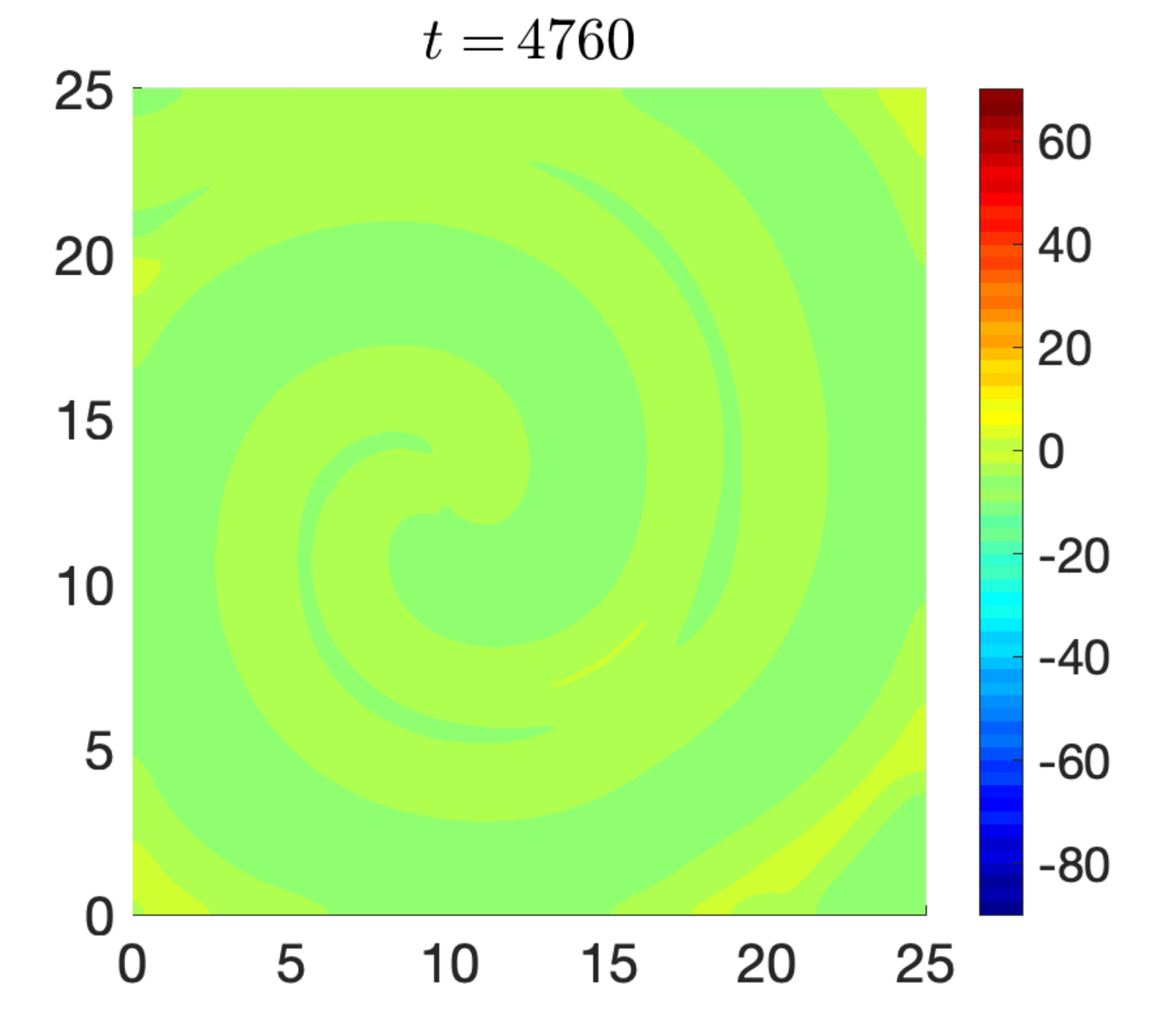}} 
\caption{A spiral wave simulated with the monodomain model \eqref{monodomain} with normal parameter settings. At $t=4s$ the parameter settings are switched to the ones used to produce Figure \ref{fig:comparison}(b), leading to cardiac death. The time $t$ is measured in milliseconds.}\label{fig:spiraldeath}
\end{figure}

Using a standard S1-S2 protocol to produce a spiral wave (see for example \cite{Vandersickel1}) with the same initial stimulus as above and with normal parameter settings, we get the spiral wave depicted in Figure \ref{fig:spiraldeath}(a). Switching the parameters to that of Figure \ref{fig:comparison}(b) after 4 seconds, which can be seen as a sudden onset of the effect of a drug, leads to cardiac death, see Figure \ref{fig:spiraldeath}(b)-(d). 

\section{Discussion and conclusion}

The main focus of this paper is to better understand the TP06 model~\cite{TP06} through bifurcation analysis. This is done by a slight reformulation of the model, which is smooth and regular enough for bifurcation theory to be applicable. To this end we reformulated the sodium current $I_\mathrm{Na}$ and fixed the intracellular potassium concentration $[K]_i$ in the TP06 ODE system. This modification enables a proper mathematical analysis and a fairly good approximation of the TP06 model~\cite{TP06} by the modified one. As demonstrated, fixing the intracellular potassium concentration $[K]_i$ by the initial intracellular potassium concentration of the TP06 model~\cite{TP06} shows similar impact on the dynamics of both models. This finding was further validated by a detailed bifurcation analysis. 

The bifurcation analysis of the modified TP06 model detected parameter values and ranges which induce EADs by a reduction of the rapid and slow potassium currents, $I_\mathrm{K_r}$ and $I_\mathrm{K_s}$, respectively. These findings indicate similar dynamics for the TP06 model~\cite{TP06}, and that this is indeed the case was demonstrated through several examples for both the single-cell and monodomain TP06 model. Thus, the bifurcation analysis of the modified model can be used to predict dynamics in the TP06 model. 

The modification does not give a perfect one-to-one relation between the two models, which is to be expected when one alters a model. It is however interesting that the modified model predicts the dynamics of the original TP06 model both for normal APs and EADs as well as it does. 

Finally, we would like to point out that the bifurcation analysis encodes the behaviour of the underlying autonomous ODE system. However, the external stimulus and the initial values play a crucial role in whether certain dynamics actually occur or not. In fact, the stimulus determines whether the trajectory reaches the basin of attraction and if dynamics like EADs appear. Examining which dynamics actually do appear in a cardiac model is a general issue and the investigation of synchronisation is highly important due to the risk that the global heart's dynamics may shift the trajectory into a dangerous region.







\vspace{0.5cm}
\noindent \small{\textbf{Conflict of interest.} The authors declare no conflict of interest.}
\\[2ex]
\small{\textbf{Acknowledgements.} 
A.E., supported in 2020 by the Kristine Bonnevie scholarship 2020 of the Faculty of Mathematics and Natural Sciences, University of Oslo, during his research stay at Lund University, wishes to thank Erik Wahl\'en and the Centre of Mathematical Sciences, Lund University, Sweden for hosting him. Since March 2021, A.E. was supported by the DFG under Germany's Excellence Strategy – MATH$^+$: The Berlin Mathematics Research Center (EXC-2046/1 – project ID: 390685689) via the project AA1-12$^*$. 
}
\bibliography{modelling_analysis}

\end{document}